\documentclass[review]{elsarticle}

\usepackage{hyperref}

\usepackage{amsfonts}
\usepackage{latexsym}
\usepackage{amsmath}

\usepackage{cases}
\usepackage{tikz}
\usepackage{graphicx}
\usepackage{graphicx}
\usepackage{epstopdf}

\newcommand{\Proof}{{\bf Proof}\hspace{5 mm}}
\newtheorem{Theorem}{Theorem}
\newtheorem{Lemma}{Lemma}

\newtheorem{Remark}{Remark}

\allowdisplaybreaks[4]

\journal{.}

\begin{document}

\begin{frontmatter}

\title{Stability and Uniqueness of Global Solutions to  Euler Equations with Exothermic Reaction \tnoteref{mytitlenote}}
\tnotetext[mytitlenote]{This paper is
supported by  National Natural Science Foundation of China (Grant No. 11426185),
Fundamental Research Funds for the Central Universities (XDJK2014C075 and SWU113062)
and  China Scholarship Council.}

\author{Kai Hu\fnref{myfootnote}}
\address{School of Mathematics and Statistics, Southwest University, Chongqing, 400715, China}
\fntext[myfootnote]{E-mail address: hukaimath@swu.edu.cn}




\begin{abstract}
We consider the Cauchy problems of a non-strictly hyperbolic system
which describes the  compressible Euler fluid with  exothermic reaction.
 In this paper a Lyapunov-type functional is constructed
for balance laws. By analysis of the flow generated by front tracking method,
we prove the well-posedness theorems
and present the local features
of  global  solutions.
\end{abstract}


\end{frontmatter}


\section{Introduction}


We are concerned about the Cauchy problem of the exothermically reacting compressible Euler equations
:
\begin{equation}\label{ZND}
  \begin{aligned}
  &v_t - u_x =0, \\
  &u_t + p_x =0 , \\
  &E_t + (pu)_x = qY\phi (T) , \\
  &Y_t= -Y\phi (T),
  \end{aligned}
\end{equation}
where $v,\ u,\ T,\ p$ and $E$ respectively represent specific volume, velocity, temperature, pressure and
total energy of compressible fluid.
This system arises from combustion theory in continuum physics.
Constant $q>0$ denotes unit binding energy released by combustion.
$Y$ represents the mass fraction
of reactant in mixed fluid, where $0\leq Y\leq 1$.
Assume that adiabatic constant $\gamma$ and specific heat constant $c$
of the mixture are independent of $Y$.
 The procedure of releaseing chemical energy is controlled by reaction-rate function $\phi(T)$, which is a
$C^1$ function w.r.t. $T$. Suppose that $\phi(T)$ is positive and increasing on $[0,+\infty)$.
For instance it is of Arrhenius form
$\phi(T)=T^\alpha \text{e}^{-\beta/T}$ with constants  $ \alpha,\beta>0 .$
In present paper, we prescribe the initial data
\begin{equation}\label{initial data}
  U_0(x)=(v_0 (x), u_0 (x), E_0 (x), Y_0 (x)),
  \ \ x \in \mathbb{R}
\end{equation}
with
$Y_0 (x) \in L^1(\mathbb{R})$ and variation $T.V.\left(U_0\right) < +\infty .$

System (\ref{ZND}) is a typical hyperbolic balance law. It can be formulated by
\begin{equation}\label{Balance Law}
  U_t + F(U)_x = G(U)
\end{equation}
where $U=(v,u,E,Y)$.
Well-posedness of solutions to balance laws was already
 established under particular circumstances.
See \cite{DafermosHsiao82,AmadoriGosseGuerra02,AmadoriGuerra02,AmadoriGuerra99,Christoforou06,ColomboGuerra07}
and the references therein.
Generally speaking, strictly hyperbolic balance laws has local weak solutions with bounded variation.
Their variation may blow up within finite time due to damping effect. Therefore,
in order to extend them to global solutions, one must impose appropriate
 dissipative conditions on source $G$.
Specifically let matrix
$$\tilde{G}=R^{-1}(\mathbf{0})\text{d}G(\mathbf{0})R(\mathbf{0}),$$
where matrix $R(U)$ consists of linearly independent right eigenvectors of $dF(U)$.
Dafermos and Hsiao\cite{DafermosHsiao82} observe that when $\tilde{G}$ is strictly diagonally dominant, system (\ref{Balance Law}) admits a global BV solution, whose variation
decreases exponentially as $t\rightarrow\infty$.
Concerning non-autonomous source $G=G(x,U)$, Amadori and collaborators\cite{AmadoriGosseGuerra02} proposed a non-resonance condition
and a smallness condition on $\omega$,
$$|G(\cdot,U)|+||\nabla_{\mathbf{u}}G(\cdot,U)|| \leq \omega \in L^1 \cap L^\infty (\mathbb{R})$$
where $||\cdot||$ denotes the Euclidian norm.
Such conditions imply global well-posedness.
Whereas our system (\ref{ZND}) is ruled out of previous situations in early research,
because it is not strictly hyperbolic.
 Two eigenvalues of Jacobian $dF(U)$ completely coincide. It is known that resonance
of different fields possibly leads to loss of uniqueness of entropy solutions.
Although one reduce (\ref{ZND}) to $3\times3$ strictly hyperbolic system
 by separating unknown $Y$ from equations (\ref{ZND}), the diagonal dominance
condition still fails to hold.
So well-posedness  of our system indeed possesses new difficulties that differ
from the previous works.

The Cauchy problem (\ref{ZND})(\ref{initial data}) has been intensively discussed in
recent research.
Main difficulty of system (\ref{ZND}) is
governing the damping induced by term $Y\phi(T)$, which arises from exothermical
reaction and amplifies any oscillation.
 To achieve this, \cite{ChenWagner03} demands  initial temperature satisfies
\begin{equation}\label{condition phi(T)}
    \phi(T_0(x))>\phi >0 ,\ \ \ \forall x \in \mathbb{R},
\end{equation}
with some constant $\phi$. Then
existence of global entropy solutions  was verified in \cite{ChenWagner03,Hu18}
if $U_0$ is of small total variation.
Related results on Navier-Stokes equations in combustion
 are referred to \cite{Chen92,ChenHoffTrivisa03,Wang04,DonatelliTrivisa06}.
However, uniqueness and stability of global solutions remain unknown.
Our paper is intended to give the complete conclusions on this problem.

In this paper, we consider the global solution $U(t,x)$ constructed  by fractional step
wave front tracking scheme, which determines a semigroup trajectory.
On stability issue, we must refer the pioneering work by Bressan, Liu
and Yang (cf. \cite{BressanLiuYang99,LiuYang99}).
They devised an elegant functional
$\Phi$ for conservation laws to analyze the $L^1$ distance of
distinct trajectories.
Unfortunately this functional is invalid for system (\ref{ZND}), since the exothermical damping
removes the rough monotonicity of $\Phi$.
It is  tough to construct an appropriate functional instead of $\Phi$.
Inspired by the argument in \cite{BressanLiuYang99},
we improve the weight included in Lyapunov functional, and then establish the
stability theorem for system (\ref{ZND}) as follows.

\begin{Theorem}\label{Theorem stability}
Assume that T.V.$(U_0)$
and $||Y_0||_{\infty}$ are sufficiently small,
moreover, condition (\ref{condition phi(T)}) holds. Then according to fractional step wave front tracking scheme,
 there exists a family of  domains $\{\mathcal{D}_t\}_{t\geq 0}$,
a unique semigroup $\mathcal{P} :\ [0,+\infty)\times\mathcal{D}_\infty
\rightarrow \mathcal{D}_\infty$ and a constant
$L>0$ such that \\
\textnormal{(i)} $U(t,x) \doteq \mathcal{P}_t (U_0)$ with $U_0  \in \mathcal{D}_0$ is an entropy weak solution to problem (\ref{ZND})({\ref{initial data}}); \\
\textnormal{(ii)} for any $U_0, V_0 \in \mathcal{D}_0$ and $ t, s\geq 0$, there holds
\begin{align}
\label{Semigroup P}   & \mathcal{P}_0(U_0)=U_0,\ \ \ \ \ \ \ \mathcal{P}_t(\mathcal{P}_s (U_0))= \mathcal{P}_{s+t}(U_0), \\
\label{Lipschitz P}   & ||\mathcal{P}_t (U_0)- \mathcal{P}_s (V_0)||_{L^1}
   \leq L(||U_0- V_0||_{L^1} + |t-s|).
\end{align}

\end{Theorem}

\begin{Remark}
We construct domain $\mathcal{D}_t$ by (\ref{Domian D}) in a general setting.
Classic stability theory regards the constant state,especially  $U\equiv\mathbf{0}$, as a equilibrium
point for balance laws. However, condition (\ref{condition phi(T)}) in this paper causes
the temperature  is globally  away from 0.
 We thereby replace the equilibrium state with a Riemann data $U_\infty$ (see Subsection \ref{Subsection Evolution in D}),
which may connect two distinct states at $x=\pm\infty$. Theorem 1 demonstrates that
the flow close to $U_\infty$ is still stable under small perturbation in $BV\bigcap L^1$ space.
\end{Remark}

Furthermore, we find that the local characteristics of solution $U(t,x)$
are completely determined by either conservation laws (\ref{conservation law})
or transport equations (\ref{transport eq}).
To illustrate this, we introduce
an entropy solution $U_\xi^C(t,x)$ to  (\ref{conservation law})
and a solution $U_\xi^T(t,x)$ to (\ref{transport eq}),
then derive the following theorem.

\begin{Theorem}
Assume that T.V.$(U_0)$
and $||Y_0||_{\infty}$ are sufficiently small,
moreover, condition (\ref{condition phi(T)}) holds.
Let $\mathcal{P}$ be a semigroup defined in Theorem 1,
and $\hat{\lambda}$ be a positive constant larger than all characteristic speeds of system (\ref{ZND}).
If $U(t,x)= \mathcal{P}_t (U_0)$ with   $U_0 \in \mathcal{D}_0$,
then for every $s\geq0$ and $a<\xi<b$, there holds
\begin{align}
\label{Local character 1}
&\mathop{\textnormal{lim}}_{\theta \rightarrow 0^+}
\frac{1}{\theta}\int_{\xi-\theta \hat{\lambda}}^{\xi+\theta \hat{\lambda}}
||U(s+\theta,x)-U_\xi^C(\theta,x)||dx =0,  \\
\label{Local character 2}
&\mathop{\overline{\textnormal{lim}}}_{\theta \rightarrow 0^+}
\frac{1}{\theta}\int_{a+\theta \hat{\lambda}}^{b-\theta \hat{\lambda}}
||U(s+\theta,x)-U_\xi^T(\theta,x)||dx  \leq
C_0 \{T.V.(U(s);(a,b))\}^2 ,
\end{align}
where  $C_0$ is a positive constant. Conversely, if mapping $U:\ [0,+\infty)
\rightarrow \mathcal{D}_\infty $  is  Lipschitz continuous in $L^1-$topology and satisfies
(\ref{Local character 1})(\ref{Local character 2})
for almost every $s\geq0$ and $a< \xi<  b$ , then $U(t)$ must coincide with the semigroup trajectory
$\mathcal{P}_t (U_0)$.
\end{Theorem}

Theorem 2 reveals that the entropy solution  possessing local features (\ref{Local character 1})(\ref{Local character 2})
is actually unique, no matter what approach we adopt.

\section{Preliminary}

The author \cite{Hu18} developed the existence theorem for general mixed Euler fluid.
First of all, we directly state the algorithm and the estimates of Glimm functional for system (\ref{ZND}).
As mentioned in Section 1, assume thermodynamics parameters $\gamma$ and $c$ remain constants
throughout the paper.

\subsection{Fractional Step Wave-Front Tracking Scheme}\label{Def FSFT scheme}

The strategy of this scheme is  separating the effect of source term from original system,
and then transforming it into a linear problem.
Now we initiate the procedure by simplifying initial data $U_0$.
For each $\varepsilon>0$, choose suitable $U^{\varepsilon}_0(x)$ as an approximation of $U_0 (x)$ such that
\begin{itemize}
  \item $U^{\varepsilon}_0(x)$ is piecewise constant, and has finite discontinuous points
        $\{x_\alpha\} \subset [-1/\varepsilon, 1/\varepsilon]$ .
  \item $T.V.\left(U^{\varepsilon}_0\right) \leq T.V.\left(U_0\right)$ .
  \item $\int_{-\frac{1}{\varepsilon}}^{\frac{1}{\varepsilon}} ||U^{\varepsilon}_0(x)-U_0(x)|| dx < \varepsilon.$
\end{itemize}

\noindent
Then divide time interval $[0, 1/\varepsilon]$ into $N$ subintervals, such that
division points
$$t_0= 0,\ \ t_1= \varepsilon,\ \ t_2= 2\varepsilon,\ \ \cdot\cdot\cdot \ \ t_N= N\varepsilon.$$
Additionally set $t_{N+1}=+\infty$.

In each strip domain
$\Omega_k= \{ \ (t,x)\ |\ t_k\leq t<t_{k+1},\ \ x \in \mathbb{R}\}$ ($k=0,1,\cdot\cdot\cdot, N)$,
 we approximately solve the initial value problem of conservation law
\begin{equation}\label{Eq: non reaction step}
  \begin{aligned}
     &U_t + F(U)_x   =0,\ \ \ \ \ \ (t,x)\in \Omega_k ,\\
     &U(t_k,x)  = U^\varepsilon (t_k,x)
   \end{aligned}
\end{equation}
by standard wave-front tracking method. See \cite{Bressanbook,Dafermosbook} for this method in detail.
Define the solution to (\ref{Eq: non reaction step}) by
$$U^\varepsilon (t,x)\doteq \mathcal{S}_t ^\varepsilon (U^\varepsilon (t_k,x)),\ \ \ \ (t,x)\in \Omega_k,$$
where $\mathcal{S}_t ^\varepsilon$ is the approximate solution operator of (\ref{Eq: non reaction step}) .
Notice that $Y(t,x)\equiv Y(t_k,x)$ for every point $(t,x)\in \Omega_k$.
When $t=t_k$, the chemical reaction will be activated.
Therefore consider an ordinary differential equations,
\begin{equation}\label{Eq2: reaction ODE}
\begin{aligned}
   &U_t= G(U), \ \ \ \ \ t \in [0,\varepsilon],\\
   &U(0,x)= U^{\varepsilon}(t_k-,x)
\end{aligned}
\end{equation}
where $U^{\varepsilon}(t_k-,x)\doteq\text{lim}_{t \rightarrow t_k^-} U^{\varepsilon}(t,x)$.
We define the linearly approximate solution by
$$ \mathcal{T}_{t} (U(0,x))  \doteq U(0,x) + G(U(0,x))t, $$
where $\mathcal{T}_{t}$ denotes the approximate solution operator of (\ref{Eq2: reaction ODE}).
Then assign the value of $U^\varepsilon (t,x)$ after reaction by
\begin{equation}\label{U(t_k,x)}
 \begin{aligned}
   U^{\varepsilon}(t_k,x) & \doteq \mathcal{T}_{\varepsilon} (U^{\varepsilon}(t_k-,x)) \\
    & = U^{\varepsilon}(t_k-,x) + G(U^{\varepsilon}(t_k-,x))\varepsilon.
 \end{aligned}
\end{equation}

In summary, we formulate the $\varepsilon-$approximate solution to system (\ref{ZND}) as follow,
\begin{equation}\label{Def P}
U^\varepsilon (t,x) = \mathcal{P}_t ^\varepsilon (U^\varepsilon_0(x))\doteq
\mathcal{S}_{t-k\varepsilon} ^\varepsilon(\mathcal{T}_{\varepsilon}\mathcal{S}_\varepsilon ^\varepsilon )^k(U^\varepsilon _0(x)),
\end{equation}
where integer $k=\text{min}\{[t/\varepsilon], N\}$.
 $\mathcal{P}_t^\varepsilon$ is actually a composition of mappings $\mathcal{T}_\varepsilon$,
$\mathcal{S}_\varepsilon^\varepsilon$, etc..
It was proved in \cite{Hu18} that as $\varepsilon\rightarrow0$,
the sequence of $U^\varepsilon (t,x)$ converges to an entropy weak solution $U(t,x)=\mathcal{P}_t  (U_0)$ of system (\ref{ZND}).


\subsection{Evolution of Approximate Solutions in $\mathcal{D}_\infty$}
\label{Subsection Evolution in D}

In this subsection we will explain how the flow of  $\mathcal{P}_t ^\varepsilon (U^\varepsilon _0(x))$
evolves as $t$ increases.
It is necessary to understand some fundamental properties of system (\ref{ZND}) in advance.
The Jacobian $dF(U)$ has four eigenvalues, i.e.
$$\lambda_1 = -\sqrt{-\gamma p_v} , \ \ \lambda_2 =0, \ \ \lambda_3 = \sqrt{-\gamma p_v}, \ \ \lambda_4 =0, $$
where $p=p(v,u,E)$. So we say that system (\ref{ZND}) is not strictly hyperbolic.
The waves of the first three fields, called elementary waves, contain shock, rarefaction wave and contact discontinuity;
the wave of the 4th field, called $Y-$wave, is of speed $\lambda_4 =0$.
To $\varepsilon-$approximate solutions, we add a family of non-physical waves (abbr. \textit{NP}) moving at fixed speed  $\hat{\lambda}$.
By construction, total strength of $NP$ waves is less than $\varepsilon$. See \cite{Hu18,Bressanbook,Dafermosbook}.

Next introduce some necessary functionals.
They  are developed from  the classic functionals in \cite{ChenWagner03,Bressanbook}.
For each $U(t,x)$, let functional
$$\mathcal{V}(U(t))\doteq \sum (|\alpha| + |\sigma| +M|\delta|), $$
where $|\alpha|$, $|\sigma|$
and $|\delta|$  respectively denote the strengths of elementary waves,  \textit{NP} wave
and $Y-$wave  in $U(t)$.
Positive constant $M$ is used to remove bad terms arising from estimates of Glimm functional;
see Lemma 8 in \cite{Hu18}.
In order to distinguish $Y-$waves which locate in different positions at time $t$,
we write $|\delta_\alpha|$ with subscript $\alpha$.
It means that this $Y-$wave front of strength $|\delta_\alpha|$ intersects with a front of strength $|\alpha|$
 at some point $(\bar{t}, \bar{x})$,  where $\bar{t}\leq t$.
We define
$|\delta_{\beta}|$, $|\delta_{\sigma}|$, $|\bar{\delta}_{\alpha}|$, $|\bar{\delta}_{\beta}|$ etc. in analogous way.
Furthermore, write
\begin{equation}\label{Def Q(U)}
\begin{aligned}
  Q(U(t))\doteq &  \mathop{\sum}\limits_{\text{App}} |\alpha|\cdot| \beta| +
 \mathop{\sum}\limits_{\text{App}}|\alpha| \cdot M|\delta _{\beta}|+
\sum M|\delta_{\alpha}|\cdot M|\delta_{\beta}|  \\
    & +  \mathop{\sum}\limits_{\text{App}} |\sigma|\cdot| \beta | +
    \mathop{\sum}\limits_{\text{App}}|\sigma| \cdot M|\delta_{\beta}|+
\sum M|\delta_{\sigma}|\cdot M|\delta_{\beta}|
\end{aligned}
\end{equation}
to measure the interaction of any two approaching waves except that of \textit{NP}
 wave and $Y-$wave. We require
the term  $M|\delta_{\alpha}|\cdot M|\delta_{\beta}|$ is included in $ Q(U(t))$
if corresponding $|\alpha |\cdot|\beta|$  exists in  $\mathop{\sum}\limits_{\text{App}} |\alpha|\cdot| \beta|$.
Analogously whether $M|\delta_{\sigma}|\cdot M|\delta_{\beta}|$ is included in $ Q(U(t))$ depends on
$|\sigma|\cdot| \beta |$.

Define modified Glimm functional
$$\mathcal{F}(U(t))\doteq \mathcal{V}(U(t))+C\cdot Q(U(t)),$$
where positive constant $C$ is so large that estimate (\ref{Boundness of F(t)}) holds.
To generate a Lipschitz continuous semigroup, we require that $Y_0 \in L^1\bigcap L^{\infty}(\mathbb{R})  $.
However, $||Y_0||_{L^1}$ may be arbitrarily large.

Suppose there exist constants $\phi$ and $\psi$ such that
$0< \phi < \phi(T_0) < \psi.$
Thus it is reasonable to assert that
\begin{equation}\label{Boundness of phi(T)}
   0< \phi \leq \phi(T(t,x)) \leq \psi
\end{equation}
holds for all $t>0$ provided $U (t,x)$ is sufficiently close to initial $U_0(x)$ in $L^\infty(\mathbb{R})$.
It was proved
that the  Glimm functional $\mathcal{F}$ designed for system (\ref{ZND}) is possibly increasing but eventually bounded, provided initial
data is suitably small. Here we directly state the following estimates of $\mathcal{F}$
in Lemma 1.
Also see conclusions (41)(42) in \cite{Hu18}.

\begin{Lemma}
Assume that $T.V.(U_0)$ and $||Y_0||_\infty$ are small enough, moreover, condition (\ref{condition phi(T)}) holds.
 Let $U^\varepsilon$
be an approximate solution defined in (\ref{Def P}), and integer $k=[t/\varepsilon]$.
Then
\begin{align}
  \mathcal{F}(U^\varepsilon(t))
    & \leq \mathcal{F}(U^\varepsilon(t_k-)) \textnormal{exp}\left\{B ||Y_0||_{\infty}
        \textnormal{e}^{-\phi k\varepsilon}\varepsilon \right\}  \nonumber\\
    & \leq \mathcal{F}(U_0) \textnormal{exp}\left\{B ||Y_0||_{\infty}
    \mathop{\sum}_{i\varepsilon\leq t}
    \textnormal{e}^{-\phi i\varepsilon}\varepsilon \right\}  \nonumber\\
   &  \leq \mathcal{F}(U_0) \textnormal{exp}\left\{\frac{B||Y_0||_{\infty}}{\phi}\right\}  \nonumber\\
\label{Boundness of F(t)}   & \leq \bar{\mathcal{F}} ,\\
\label{Decay of Y(t)}
  Y(t) &\leq Y(0)\textnormal{e}^{ -\phi k\varepsilon }
\end{align}

\noindent for small $\varepsilon>0$, where $i\in \mathbb{N}_+$, both $B$ and $\bar{\mathcal{F}}$
are constants .
\end{Lemma}

Comparing (\ref{Boundness of F(t)}) with original result in \cite{Hu18} ,
 we only maintain the first order terms of $\varepsilon$,
because all higher order ones can be absorbed if $\varepsilon$ is small.
Since $\mathcal{F}(U)$ is equivalent to $T.V.(U)$, there exists a constant
$B^*>0$ such that
$$\mathcal{F}(U_0) < \epsilon  \text{ \ \ \ implies \ \ \ }
||Y_0||_{\infty}<B^*\epsilon$$
for each $\epsilon>0$.
Given $\epsilon$,  define
\begin{align*}
   & \epsilon(t)\doteq\epsilon\cdot\text{exp}\left\{ BB^* \epsilon
    \int_{0}^{t}
    \text{e}^{-\phi s}\text{d}s \right\}
    =\epsilon\cdot\text{exp}\left\{ \frac{BB^*\epsilon}{\phi}
    (1-\text{e}^{-\phi t}) \right\},   \\
   &\epsilon_\infty \doteq\epsilon\cdot\text{exp}\left\{\frac{BB^*\epsilon}{\phi} \right\}.
\end{align*}

\noindent  Set vector $\mathcal{U} \doteq (v,u,E)$.
State $U_{\infty}=(\mathcal{U}_\infty(x),0)$ denotes a Riemann initial data with two pieces of constant states, i.e.
\begin{equation*}
    \mathcal{U}_\infty(x)=\left\{
    \begin{aligned}
        &\mathcal{U}_{r}, \ \ \ x>0, \\
        &\mathcal{U}_{l}, \ \ \ x<0.
    \end{aligned}\right.
\end{equation*}
Evidentally $G(U_{\infty})=\mathbf{0}$ in (\ref{ZND}) and (\ref{Balance Law}).
Now we can explicitly define the domains where
the flow $\mathcal{P}_t^\varepsilon(U)$ evolves.
In fact, for every  $t\geq0$, set
\begin{align*}
  &\mathcal{M}_t\doteq  \textbf{cl} \left\{ \bar{U} \in L^1(\mathbb{R};\mathbb{R}^4)\ \right. |\
\bar{U} \text{ is piecewise constant}, \\
&\ \ \ \ \ \ \ \ \ \ \ \ \ \ \ \ \ \ \ \ \ \ \ \
 \left. \ \mathcal{F}(\bar{U}+U_{\infty}) < \epsilon (t) ,
 \ ||\bar{Y}||_{\infty}<B^*\epsilon \text{e}^{-\phi t}\right\}, \\
  &\mathcal{M}_\infty\doteq  \textbf{cl} \left\{ \bar{U} \in L^1(\mathbb{R};\mathbb{R}^4)\ \right. |\
\bar{U} \text{ is piecewise constant}, \\
&\ \ \ \ \ \ \ \ \ \ \ \ \ \ \ \ \ \ \ \ \ \ \ \
 \left. \ \mathcal{F}(\bar{U}+U_{\infty}) < \epsilon _\infty ,
 \ ||\bar{Y}||_{\infty}=0 \right\},
\end{align*}
where \textbf{cl} denotes closure in $L^1$ topology. Then define domain

\begin{equation}\label{Domian D}
    \mathcal{D}_t\doteq \left\{ U\ \left|\ U= \bar{U}+U_{\infty},\ \bar{U}
\in \mathop{\bigcup}_{s\leq t}\mathcal{M}_s
\right. \right\}.
\end{equation}

\noindent These domains obviously satisfy
$\mathcal{D}_s\subset  \mathcal{D}_t\subset \mathcal{D}_{\infty}$
for $0\leq s<t$.

We complete the section by the dynamics of $\mathcal{P}^\varepsilon_t (U)$ in $\mathcal{D}_{\infty}$.

\begin{Lemma}
Assume that $\epsilon$ is sufficiently small.
If $U \in \mathcal{D}_{s}$ is a piecewise constant function, then there exists a
 positive $\bar{t}$, which depends on $U, \epsilon$ and $ s$, satisfying that \\
\textnormal{(1)} $\mathcal{T}_t (U) \in \mathcal{D}_{s+t}$ for $t \in (0, \bar{t} )$;\\
\textnormal{(2)} $\mathcal{S}^\varepsilon_t (U) \in \mathcal{D}_{s}$,
$\mathcal{P}^\varepsilon_t (U) \in \mathcal{D}_{s+t}$ for any $t\in (0,+\infty)$ , provided $\varepsilon \ll \bar{t}$.
\end{Lemma}

{\Proof}
1. If $U \in \mathcal{D}_{s}$,  there exists
some $r<s$ such that $U \in \mathcal{M}_{r}$. This gives
\begin{align*}
   & \mathcal{F}(U) < \epsilon (r)=\epsilon\cdot\text{exp}\left\{ \frac{BB^*\epsilon}{\phi}
    (1-\text{e}^{-\phi r}) \right\}, \\
    & ||Y||_{\infty}<B^*\epsilon \text{e}^{-\phi r}.
\end{align*}

\noindent Thanks to  Lemma 10 in \cite{Hu18} and the definition of $\mathcal{T}_t$,
 we have the following estimates.
\begin{equation}
\begin{aligned}
&  \mathcal{F}(\mathcal{T}_t(U))
     \leq \mathcal{F}(U) \text{exp}\left\{B ||Y||_{\infty}
        \text{e}^{-\phi t}t \right\}\\
&  \ \ \ \ \ \ \  \ \ \ \   \leq \mathcal{F}(U) \text{exp}\left\{B ||Y||_{\infty}t \right\}, \\
&  \mathcal{T}_t Y =Y-tY\phi(T) \leq Y\text{e}^{ -\phi t },
\end{aligned}
\end{equation}
where $\mathcal{T}_t Y$ denotes the $Y-$component of $\mathcal{T}_t(U)$.
Immediately, they imply
\begin{align}
\label{Estimate 2}
&  \mathcal{F}(\mathcal{T}_t(U))
     \leq \epsilon\cdot\text{exp}\left\{ \frac{BB^*\epsilon}{\phi}
    (1-\text{e}^{-\phi r})+ B ||Y||_{\infty}t  \right\}, \\
\label{Estimate 5}
&  ||\mathcal{T}_t Y||_{\infty} \leq ||Y||_{\infty}\text{e}^{ -\phi t }
<B^*\epsilon \text{e}^{-\phi (r+t)} \ \  \ \ \text{if}\ \ t\phi<1.
\end{align}

\noindent Set function
$$\varphi(t)=\frac{\phi t}{1- \text{e}^{ -\phi t }}.$$
Since $\varphi(0)\rightarrow 1$ as $t\rightarrow 0$, we select a positive $\bar{t}$
such that
$$||Y||_{\infty} \varphi(t) < B^*\epsilon \text{e}^{-\phi r},$$
$$\text{i.e.}\ \ \ \ B||Y||_{\infty} t <\frac{BB^*\epsilon}{\phi} \text{e}^{-\phi r}
(1-\text{e}^{-\phi t}) $$
for every $t \in (0,\bar{t})$.
Thus it follows from (\ref{Estimate 2}) that
\begin{equation}\label{Estimate 3}
\mathcal{F}(\mathcal{T}_t(U)) < \epsilon (r)\text{exp}\left\{B ||Y||_{\infty}
        t\right\}
        <\epsilon\cdot\text{exp}\left\{ \frac{BB^*\epsilon}{\phi}
    (1-\text{e}^{-\phi (r+t)}) \right\}.
\end{equation}
(\ref{Estimate 5}) and (\ref{Estimate 3}) guarantee  $\mathcal{T}_t(U) \in \mathcal{M}_{r+t}\subset \mathcal{D}_{s+t}$.

2. Since $0<\varepsilon\ll \bar{t}$, it is deduced from standard argument of front tracking scheme that
the Glimm functional $\mathcal{F}(\mathcal{S}_t^\varepsilon (U)) $ is decreasing in $t$,
provided $T.V.(U) $
is small enough. Therefore, if $U \in \mathcal{D}_{s}$, the fact
$\mathcal{F}(\mathcal{S}_t^\varepsilon (U)) \leq \mathcal{F}(U)$
and $Y(t)\equiv Y(0)$ yields $\mathcal{S}_t^\varepsilon (U) \in \mathcal{D}_{s}$ for any $t>0$.
Recall that $\mathcal{P}_t^\varepsilon (U)= \mathcal{S}_{t-k\varepsilon}^\varepsilon
(\mathcal{T}_\varepsilon  \mathcal{S}_\varepsilon^\varepsilon)^k (U)$
where $k=[t/\varepsilon]$. The preceding results of $\mathcal{T}_t$ and $\mathcal{S}_t$
directly yield $\mathcal{P}^\varepsilon_t (U) \in \mathcal{D}_{s+t}$.
\hfill{$\square$}

\section{Stability of Entropy Solution}

We will construct a functional in this section to measure the distance of two approximate solutions.
This functional is almost decreasing in time $t$, which implies the  continuous dependence of entropy
solutions on initial data.

\subsection{Lyapunov Functional}

Assume initial data $U_0$ and $ V_0$ belong to $\mathcal{D}_0$.
$U_0^\varepsilon$ and $ V_0^\varepsilon$ are corresponding $\varepsilon-$approximations.
Let
$$\mathcal{P}_t^\varepsilon(U_0^\varepsilon) = U^\varepsilon(t,x)=(\mathcal{U}_1, Y_1),$$
$$\mathcal{P}_t^\varepsilon(V_0^\varepsilon) = V^\varepsilon(t,x)=(\mathcal{U}_2, Y_2).$$
Here every state $U$ (or $V$) is decomposed into components $\mathcal{U} \doteq (v,u,E)$ and $Y$.

Consider equations (\ref{ZND}) in phase space.
The $i-$th family of Hugoniot curves is determined by the mapping $H_i$.
For given states $\mathcal{U}_1$ and $\mathcal{U}_2$, there exists unique vector
$\mathbf{q}= (q_1, q_2, q_3)$ such that
\begin{equation}\label{Def q(x)}
    \mathcal{U}_2 = H(\mathbf{q})(\mathcal{U}_1)
\doteq H_3(q_3)\circ H_2(q_2)\circ H_1(q_1)(\mathcal{U}_1)
\end{equation}
at every point $(t,x)$.
It is known that total strength $\mathop{\sum}_{i=1}^3 |q_i|$ is equivalent to $||\mathcal{U}_1-\mathcal{U}_2||$.

We introduce the  Lyapunov functional
\begin{equation}\label{Lyapunov functional}
  \Phi (U^\varepsilon(t), V^\varepsilon(t)) \doteq \mathop{\sum}\limits_{i=1}^{3}
        \int_{-\infty}^{+\infty} \left(|q_i(t,x)|+ \kappa_1 |Y_1(t,x) - Y_2(t,x)| \right)W_i(t,x) dx.
\end{equation}

\noindent This functional is developed from the one in \cite{BressanLiuYang99}. It is worth pointing out that,
in contrast with $\Phi$ in  \cite{BressanLiuYang99}, we add the term of $Y$ and modify the weight $W_i$
in order to overcome the damping effect.
Specifically, let $\alpha$ be a jump of the $k_\alpha-$th  family ($1\leq k_\alpha\leq 3$),
whose location is $(t,x_\alpha)$  and strength $|\alpha|$ .
 Let $\delta$ be a jump of $Y-$wave, whose  strength is $|\delta|$.
Notation $\mathcal{J}(U)$ denotes a set of all jumps of $U$ for given time $t$.
Set $\mathcal{J}\doteq\mathcal{J}(U^\varepsilon)\bigcup \mathcal{J}(V^\varepsilon)$.
Then the weight
$$
W_i(t,x)\doteq 1+ \kappa_2 A_i(t,x) + \kappa_3  \{Q(U^\varepsilon)+Q(V^\varepsilon)\} +\kappa_4 B(t),
$$

\noindent where for $i=1,3,$
\begin{align*}
  A_i(t,x)= &  \left[\sum\limits_{
                         \begin{subarray}{c}
                         x_\alpha <x,\ i<k_\alpha \leq 3, \\
                         \alpha \in \mathcal{J}
                          \end{subarray}
                          }
                +\sum\limits_{
                        \begin{subarray}{c}x_\alpha >x,\ 1\leq k_\alpha <i, \\
                        \alpha \in \mathcal{J}
                          \end{subarray}
                          }
                \right] |\alpha|
+  \sum\limits_{\delta \in \mathcal{J}} M|\delta|\\
& + \left\{
    \begin{aligned}
       & \left[\sum\limits_{
            \begin{subarray}{c}
                  k_\alpha =i ,  x_\alpha <x,\\
               \alpha \in \mathcal{J}(U^\varepsilon)\\
            \end{subarray}
            }
       + \sum\limits_{
            \begin{subarray}{c}
              k_\alpha =i ,  x_\alpha >x \\
                \alpha \in \mathcal{J}(V^\varepsilon)
            \end{subarray}
            }
      \right] |\alpha|
      \ \ \ & \mbox{if}\ q_i(x)<0 ,\\
      & \left[\sum\limits_{
          \begin{subarray}{c}
                  k_\alpha =i ,  x_\alpha <x,\\
               \alpha \in \mathcal{J}(V^\varepsilon)\\
            \end{subarray}
            }
       + \sum\limits_{
           \begin{subarray}{c}
              k_\alpha =i ,  x_\alpha >x \\
                \alpha \in \mathcal{J}(U^\varepsilon)
            \end{subarray}
            }
      \right] |\alpha|
      & \mbox{if}\ q_i(x) >0,
    \end{aligned}
    \right. \\
  A_2(t,x)= &  \left[\sum\limits_{
                         \begin{subarray}{c}
                         x_\alpha <x,\ k_\alpha = 3, \\
                         \alpha \in \mathcal{J}
                          \end{subarray}
                          }
                +\sum\limits_{
                        \begin{subarray}{c}x_\alpha >x,\  k_\alpha=1 , \\
                        \alpha \in \mathcal{J}
                          \end{subarray}
                          }
                \right] |\alpha|
+  \sum\limits_{\delta \in \mathcal{J}} M|\delta|,\\
B(t)=& \mathop{\sum}\limits_{j\varepsilon > t}
        (||Y_1 (0,\cdot)||_{\infty}+||Y_2 (0,\cdot)||_{\infty}) \text{e}^{-\phi j\varepsilon}\varepsilon,
\end{align*}

\noindent and $Q(U)$ defined in (\ref{Def Q(U)}). The constants $\kappa_i \ ( 1\leq i \leq 4) $ will be specified in the sequel.
Especially $\kappa_3 = C \cdot\kappa_2$. Recall that $C$ is a large constant in the definition of $\mathcal{F}(U)$.

We claim  $1\leq W_i \leq 2$ if $\epsilon$ is small enough.
By the boundness of $W_i$,
one can verify that $L^1$ distance of $U^\varepsilon(t)$ and $V^\varepsilon(t)$ is equivalent to $ \Phi (U^\varepsilon(t), V^\varepsilon(t))$,
i.e.
\begin{equation}\label{equivalence metric}
    \frac{1}{C_1}||U^\varepsilon - V^\varepsilon||_{L^1} \leq  \Phi (U^\varepsilon, V^\varepsilon) \leq C_1||U^\varepsilon - V^\varepsilon||_{L^1}
\end{equation}
for some constant $C_1>1$.
To understand the property of $\Phi$, let us prove the monotonicity of $W_i$ at first.

\begin{Lemma}\label{Lemma Wi decreasing}
For every $x \in \mathbb{R}$, $W_i(x,t)$ is piecewise constant and decreasing in $t$.
\end{Lemma}

{\Proof}
When $t \in (t_k, t_{k+1})$, all the terms w.r.t. $Y$ are unchanged. Thus $W_i$ is
equivalent to the weight $W_i(x)$ defined for conservation laws in \cite{BressanLiuYang99}.
This is certainly decreasing in $t$. It suffices to verify the case  $t=t_k$.
Set
$$W_i^+\doteq W_i(t_k,x),\ \ \ \ W_i^-\doteq W_i(t_k-,x).$$
Using (\ref{Boundness of F(t)}), we calculate that
\begin{equation}\label{Estimate W+ - W-}
\begin{aligned}
  W_i^+ - W_i^- & \leq  \kappa_2 \bar{\mathcal{F}}B (||Y_1 (0,\cdot)||_{\infty} + ||Y_2 (0,\cdot)||_{\infty})  \cdot  \text{e}^{-\phi k\varepsilon}\varepsilon  \\
   &\ \ \  - \kappa_4  (||Y_1 (0,\cdot)||_{\infty} + ||Y_2 (0,\cdot)||_{\infty})
       \cdot  \text{e}^{-\phi k\varepsilon} \varepsilon \\
   & \leq - \frac{1}{2}\kappa_4 (||Y_1 (0,\cdot)||_{\infty} + ||Y_2 (0,\cdot)||_{\infty})
    \cdot  \text{e}^{-\phi k\varepsilon} \varepsilon \\
   &<0
\end{aligned}
\end{equation}
if $\kappa_4 $ is so large that $\kappa_2 \bar{\mathcal{F}}B < \frac{1}{4} \kappa_4$.
The proof is completed. \hfill{$\square$}

The modified weight $W_i$ is  absolutely a key factor in the whole paper. Using its
monotonicity, we successfully prove in Lemma \ref{Lemma phi almost decreasing}
that $\Phi$ is almost deceasing. This fact implies the
stability and uniqueness of entropy  solution.

According to the scheme stated in
Subsection \ref{Def FSFT scheme}, time $t$ can be categorized into three cases as follows.
The time when approaching waves interact is denoted by $t^*$;
The time when chemical reaction is activated is denoted by $t_k$ ;
The other time is denoted by $t^\circ$ .
Obviously Lyapunov functional $ \Phi(U^\varepsilon(t),V^\varepsilon(t))$ is decreasing at time $t^*$,
because the main factor
$$\int_{-\infty}^{+\infty} \left(|q_i(t,x)|+ \kappa_1 |Y_1(t,x) - Y_2(t,x)| \right)dx$$
 in definition (\ref{Lyapunov functional})
is  continuous  at $t^*$, meanwhile, $W_i$ is decreasing steeply.
Hence it suffices to discuss the remaining cases $t=t_k$ or $t^\circ$ .

\subsection{Estimate of $\Phi(t)$ at $t_k$}

Define  the components of $U^\varepsilon$
on two banks of discontinuity $t=t_k$ by
\begin{align*}
    & \mathcal{U}_1^+ \doteq \mathcal{U}_1 (t_k , x),  \ \ \ \  \mathcal{U}_1^- \doteq \mathcal{U}_1 (t_k -, x),\\
    & Y_1^+ \doteq  Y_1(t_k , x),   \ \ \ \   Y_1^- \doteq Y_1 (t_k -, x).
\end{align*}

\noindent And $\mathcal{U}_2^\pm,  Y_2^\pm $ represent the similar quantities of $V^\varepsilon$.
Accordingly define $\mathbf{q}^+$ and $\mathbf{q}^-$ by expression
(\ref{Def q(x)}).

\begin{Lemma}\label{Lemma Estimate q_i}
At every point $(t_k,x)$, assume $\mathbf{q}^+$ and $\mathbf{q}^-$ satisfy
\begin{equation}\label{Def q+ q-}
    \mathcal{U}_2^- = H(\mathbf{q}^-)(\mathcal{U}_1^-),\ \ \
\mathcal{U}_2^+ = H(\mathbf{q}^+)(\mathcal{U}_1^+).
\end{equation}
 Then
\begin{equation}\label{Estimate q_i}
  q_i^+ - q_i ^- = \mathcal{O}(1)
  \left\{|Y_1^-  - Y_2^-|\varepsilon +
   \mathop{\sum}_{j=1}^3 |q_j^-| (Y_1^- + Y_2^-)\varepsilon\right\} .
\end{equation}
holds for $i=1,2,3$.
\end{Lemma}
The Landau symbol $\mathcal{O}(1)$ denotes a quantity whose absolute value is bounded by a constant dependent
of System (\ref{ZND}) and the domain $\mathcal{D_\infty}$.

{\Proof}
Set the differences
\begin{align*}
  &\Delta_1  \doteq Y^+_1 - Y^- _1= -Y^- _1\phi(T^-_1)\varepsilon,\\
  &\Delta_2  \doteq Y^+_2 - Y^-_2 = -Y^-_2 \phi(T^-_2)\varepsilon,
\end{align*}

\noindent and vector $\mathbf{w} \doteq (0,0,-q).$
Then it follows that
\begin{equation}\label{Def Ua+ Ub+}
\mathcal{U}_1^+= \mathcal{U}_1 ^- + \mathbf{w}\cdot \Delta_1,\ \ \
\mathcal{U}_2^+= \mathcal{U}_2^- + \mathbf{w}\cdot \Delta_2 .
\end{equation}

\noindent Equations (\ref{Def q+ q-})
and (\ref{Def Ua+ Ub+}) determine an implicit $C^2$ function
$$\mathbf{q}^+ \doteq \mathbf{q}(\mathbf{q}^-, \Delta_1,\Delta_2).$$

\noindent Notice that $\mathbf{q}(\mathbf{0},0,0)= \mathbf{q}(\mathbf{0},\Delta_1,\Delta_1)=\mathbf{0}$, and $\mathbf{q}(\mathbf{q}^-,0,0)=\mathbf{q}^-$.
Using Proposition 1 in \cite{Hu18} and Taylor's formula, we figure out
\begin{equation}
\begin{aligned}\label{Estimate q+}
   \mathbf{q}^+
   &= \mathbf{q}(\mathbf{0},\Delta_1,\Delta_2) +\mathbf{q}(\mathbf{q}^-,0,\Delta_2)
        -\mathbf{q}(\mathbf{0},0,\Delta_2) +\mathcal{O}(1)||\mathbf{q}^-||\cdot|\Delta_1|\cdot\mathbf{I}\\
   &= \mathbf{q}(\mathbf{0},\Delta_1,\Delta_2) + \mathbf{q}^-
         +\mathcal{O}(1)||\mathbf{q}^-||\cdot(|\Delta_1|+|\Delta_2|)\cdot\mathbf{I}\\
   &=  \mathbf{q}^-  +\mathcal{O}(1)|\Delta_1 - \Delta_2|\cdot\mathbf{I}
         +\mathcal{O}(1)||\mathbf{q}^-||\cdot(|\Delta_1|+|\Delta_2|)\cdot\mathbf{I}
\end{aligned}
\end{equation}
where $\mathbf{I}=(1,1,1)^T$ is a $3\times1$ vector.
By Lipschitz continuity of $\phi(T)$, we have
\begin{equation}\label{Estimate phi}
    |\phi(T_1^-)-\phi(T_2^-)| \leq
\mathcal{O}(1) ||\mathcal{U}_1^-  -  \mathcal{U}_2^-||
\leq \mathcal{O}(1) \mathop{\sum}_{j=1}^3 |q_j^-|.
\end{equation}

\noindent It follows from condition (\ref{Boundness of phi(T)}) and (\ref{Estimate phi}) that
\begin{align}
    |\Delta_1 - \Delta_2|
    & \leq |Y_1^-  - Y_2^-|\phi(T^-_1)\varepsilon +Y_2^-|\phi(T^-_1)-\phi(T^-_2)|\varepsilon \nonumber \\
\label{Estimate deltaa-deltab}
    & \leq \mathcal{O}(1)|Y_1^-  - Y_2^-| \varepsilon +\mathcal{O}(1) Y_2^-||\mathcal{U}^-_1 -\mathcal{U}^-_2||\varepsilon \\
    &\leq \mathcal{O}(1)\left\{|Y_1^-  - Y_2^-|\varepsilon +
        \mathop{\sum}_{j=1}^3 |q_j^-|  Y_2^- \varepsilon\right\} , \nonumber\\
\label{Estimate deltaa+deltab}
   |\Delta_1 |+| \Delta_2|
     &\leq   \mathcal{O}(1) (Y_1^- + Y_2^-)\varepsilon.
\end{align}

\noindent Finally (\ref{Estimate q+})(\ref{Estimate deltaa-deltab})
and (\ref{Estimate deltaa+deltab}) yield estimate (\ref{Estimate q_i}).
\hfill{$\square$}

\begin{Lemma}\label{Lemma Estimate Y_1-Y_2}
At every point $(t_k,x)$, components $Y$ of $U^\varepsilon$ and $V^\varepsilon$ satisfy that
\begin{equation}\label{Estimate Y_1-Y_2}
  |Y_1^+ - Y_2^+| - |Y_1^- - Y_2^-| \leq
  -|Y_1^-  - Y_2^-|\phi \varepsilon + \mathcal{O}(1)  \mathop{\sum}_{j=1}^3 |q_j^-|  Y_2^-\varepsilon.
\end{equation}
\end{Lemma}

{\Proof}
Note that $\phi(T)$ is bounded below by $\phi$. It is easy to check that
\begin{align*}
  |Y_1^+ - Y_2^+| = &|Y_1^-(1-\phi(T_1^-)\varepsilon) - Y_2^-(1-\phi(T_2^-)\varepsilon)| \\
  \leq & |Y_1^- - Y_2^-|\cdot|1-\phi(T_1^-)\varepsilon| + Y_2^-|\phi(T_1^-)-\phi(T_2^-)|\varepsilon \\
  \leq & |Y_1^- - Y_2^-|\cdot(1-\phi\varepsilon) + \mathcal{O}(1)
   \mathop{\sum}_{j=1}^3 |q_j^-|Y_2^-\varepsilon
\end{align*}
for small $\varepsilon$. This implies the estimate (\ref{Estimate Y_1-Y_2}).
\hfill{$\square$}

\begin{Lemma}\label{Lemma phi decreasing}
At every time $t_k=k\varepsilon$, Lyapunov functional is decreasing, i.e.
$$ \Phi(U^\varepsilon(t_k ),V^\varepsilon(t_k )) \leq
\Phi(U^\varepsilon(t_k -),V^\varepsilon(t_k -)).$$

\end{Lemma}

{\Proof}
By Lemmas \ref{Lemma Estimate q_i}-\ref{Lemma Estimate Y_1-Y_2} and
estimate (\ref{Boundness of F(t)})
we have
\begin{equation}\label{Estimate 1}
\begin{aligned}
   & (|q_i^+| +\kappa_1 |Y_1^+ - Y_2^+|)- (|q_i ^-|  + \kappa_1|Y_1^- - Y_2^-|)\\
  = & (|q_i^+| - |q_i ^-|) + \kappa_1 (|Y_1^+ - Y_2^+| - |Y_1^- - Y_2^-|)\\
  \leq & (-\kappa_1 \phi + \mathcal{O}(1)) |Y_1^- - Y_2^-|\varepsilon + \mathcal{O}(1)  \mathop{\sum}_{j=1}^3 |q_j^-| (Y_1^- + Y_2^-)\varepsilon \\
  \leq & - \frac{1}{2}\kappa_1 \phi  |Y_1^- - Y_2^-|\varepsilon +\mathcal{O}(1)  \mathop{\sum}_{j=1}^3 |q_j^-| (Y_1^- + Y_2^-)\varepsilon \\
  \leq &  \mathcal{O}(1)  \mathop{\sum}_{j=1}^3 |q_j^-| \cdot
   (||Y_1 (0,\cdot)||_{\infty} + ||Y_2 (0,\cdot)||_{\infty})
   \text{e}^{-\phi k\varepsilon}\varepsilon
\end{aligned}
\end{equation}
provided $\kappa_1$ is large enough.
Then  estimates (\ref{Estimate 1}) and (\ref{Estimate W+ - W-})
give
\begin{align*}
  &  \ \ \ (|q_i^+|+ \kappa_1 |Y_1^+ - Y_2^+|)W_i^+ - (|q_i^-|+ \kappa_1  |Y_1^- - Y_2^-|)W_i^- \\
  & = \left\{ (|q_i^+| - |q_i ^-|) + \kappa_1 (|Y_1^+ - Y_2^+| - |Y_1^- - Y_2^-|)\right\}W_i^+  \\
  &   \ \ \ + (|q_i^-|+ \kappa_1 |Y_1^- - Y_2^-|) (W_i^+ -W_i^-)\\
   & \leq \mathcal{O}(1) \mathop{\sum}\limits_{j=1}^3 |q_j^-| \cdot (||Y_1 (0,\cdot)||_{\infty} + ||Y_2 (0,\cdot)||_{\infty}) \text{e}^{-\phi k\varepsilon}\varepsilon \cdot W_i^+  \\
  & \ \ \ - \frac{1}{2}\kappa_4 (|q_i^-|+ \kappa_1 |Y_1^- - Y_2^-|)(||Y_1 (0,\cdot)||_{\infty} + ||Y_2 (0,\cdot)||_{\infty})\cdot  \text{e}^{-\phi k\varepsilon}  \varepsilon \\
   & \leq \mathcal{O}(1) \mathop{\sum}\limits_{j=1}^3 |q_j^-| \cdot (||Y_1 (0,\cdot)||_{\infty} + ||Y_2 (0,\cdot)||_{\infty}) \text{e}^{-\phi k\varepsilon}\varepsilon   \\
  & \ \ \   - \frac{1}{2}\kappa_4 |q_i^-|\cdot (||Y_1 (0,\cdot)||_{\infty} + ||Y_2 (0,\cdot)||_{\infty}) \cdot  \text{e}^{-\phi k\varepsilon}  \varepsilon  .
\end{align*}

\noindent Furthermore
\begin{align*}
  &  \ \ \ \mathop{\sum}\limits_{i=1}^3 (|q_i^+|+ \kappa_1 |Y_1^+ - Y_2^+|)W_i^+ -
  \mathop{\sum}\limits_{i=1}^3 (|q_i^-|+ \kappa_1  |Y_1^- - Y_2^-|)W_i^- \\
  & \leq \mathcal{O}(1) \mathop{\sum}\limits_{j=1}^3 |q_j^-| \cdot (||Y_1 (0,\cdot)||_{\infty} + ||Y_2 (0,\cdot)||_{\infty})\cdot \text{e}^{-\phi k\varepsilon}\varepsilon  \\
  & \ \ \   - \frac{1}{2}\kappa_4 \mathop{\sum}\limits_{i=1}^3 |q_i^-|\cdot
  (||Y_1 (0,\cdot)||_{\infty} + ||Y_2 (0,\cdot)||_{\infty}) \cdot
   \text{e}^{-\phi k\varepsilon}\varepsilon \\
  & \leq 0
\end{align*}
provided $\kappa_4$ is large enough.
Consequently
\begin{align*}
   & \Phi(U^\varepsilon(t_k ),V^\varepsilon(t_k )) -  \Phi(U^\varepsilon(t_k -),V^\varepsilon(t_k -)) \\
  = & \int_{-\infty}^{+\infty} \left\{
  \mathop{\sum}\limits_{i=1}^3 (|q_i^+|+ \kappa_1 |Y_1^+ - Y_2^+|)W_i^+ -
  \mathop{\sum}\limits_{i=1}^3 (|q_i^-|+ \kappa_1  |Y_1^- - Y_2^-|)W_i^- \right\} dx \\
  \leq & 0  .
\end{align*}
\hfill{$\square$}

\subsection{Estimate of $\Phi(t)$ at $t^\circ$}

At any time $ t^\circ$, jump $\alpha \in \mathcal{J}$
locates at $x_\alpha$ and moves at the speed of $\dot{x}_\alpha$.
We use the notations
\begin{align*}
    & W_i^{\alpha+} \doteq  W_i(t^\circ,x_\alpha+),&\   &W_i^{\alpha-} \doteq  W_i(t^\circ,x_\alpha-),  \\
    & q_i^{\alpha+} \doteq  q_i(t^\circ,x_\alpha+),&  &q_i^{\alpha-} \doteq  q_i(t^\circ,x_\alpha-),
\end{align*}
etc. to represent the quantities on two banks of $x_\alpha$.
Then we find $\Phi$ is almost decreasing at $t^\circ$.

\begin{Lemma}\label{Lemma phi almost decreasing}
At every $t=t^\circ$, Lyapunov functional satisfies
$$\frac{d}{dt}\Phi(U^\varepsilon(t),V^\varepsilon(t))
\leq \mathcal{O}(1) \varepsilon.
$$
\end{Lemma}

{\Proof}
We can calculate the derivative of $\Phi$ by means of the technique in \cite{BressanLiuYang99}.
Set $w_0=\mathcal{U}_1, \ w_i=H_i(q_i)\circ\cdots\circ H_1(q_1)( \mathcal{U}_1), \ w_3=\mathcal{U}_2.$
Then define $\lambda_i \doteq \lambda_i(w_{i-1},w_i)$ as the speed of the $i-$shock connecting states
$w_{i-1}$ and $w_i$. We derive that
\begin{equation*}\label{dt_Psi}
\begin{aligned}
    \frac{d}{dt}\Phi(U^\varepsilon(t),V^\varepsilon(t))
    =&\mathop{\sum}_{\alpha \in \mathcal{J}}\mathop{\sum}_{i=1}^{3}
\left\{\left(|q_i^{\alpha-}|+\kappa_1|Y_1-Y_2|^{\alpha-}\right)W_i ^{\alpha-}\right.\\
&\ \ \ \ \ \ -\left.\left(|q_i^{\alpha+}|+\kappa_1|Y_1-Y_2|^{\alpha+}\right)
    W_i ^{\alpha+}\right\}\cdot \dot{x}_\alpha \\
=&\mathop{\sum}_{\alpha \in \mathcal{J}}\mathop{\sum}_{i=1}^{3}
 \left\{\left(|q_i^{\alpha+}|+\kappa_1|Y_1-Y_2|^{\alpha+}\right)W_i ^{\alpha+}(\lambda_i^{\alpha+}-\dot{x}_\alpha )\right.\\
&\ \ \ \ \ \ -\left.\left(|q_i^{\alpha-}|+\kappa_1|Y_1-Y_2|^{\alpha-}\right)
    W_i ^{\alpha-}(\lambda_i^{\alpha-}-\dot{x}_\alpha )\right\}\\
=& \mathop{\sum}_{\alpha \in \mathcal{J}}\mathop{\sum}_{i=1}^{3} E_{i,\alpha}
\end{aligned}
\end{equation*}
where
\begin{equation*}\label{}
\begin{aligned}
E_{i,\alpha}
\doteq & \left(|q_i^{\alpha+}|+\kappa_1|Y_1-Y_2|^{\alpha+}\right)W_i ^{\alpha+}(\lambda_i^{\alpha+}-\dot{x}_\alpha )\\
&-\left(|q_i^{\alpha-}|+\kappa_1|Y_1-Y_2|^{\alpha-}\right)
W_i ^{\alpha-}(\lambda_i^{\alpha-}-\dot{x}_\alpha ).
\end{aligned}
\end{equation*}

\noindent Observe $|Y_1-Y_2|^{\alpha-}=|Y_1-Y_2|^{\alpha+}$ at $t^\circ$.
One compute that
\begin{equation*}\label{}
\begin{aligned}
E_{i,\alpha}
=& \left(|q_i^{\alpha+}|-|q_i^{\alpha-}\right)W_i ^{\alpha+}(\lambda_i^{\alpha+}-\dot{x}_\alpha )\\
&+ \left(|q_i^{\alpha-}|+\kappa_1|Y_1-Y_2|^{\alpha-}\right)(W_i ^{\alpha+}-W_i ^{\alpha-})(\lambda_i^{\alpha+}-\dot{x}_\alpha) \\
&+ \left(|q_i^{\alpha-}|+\kappa_1|Y_1-Y_2|^{\alpha-}\right)W_i ^{\alpha-}(\lambda_i^{\alpha+}-\lambda_i^{\alpha-}). \\
\end{aligned}
\end{equation*}

\noindent Let $S, R$ and $NP$ respectively denote the front of shock, rarefaction and non-physical wave.
Since $|q_i^{\alpha-}|+\kappa_1|Y_1-Y_2|^{\alpha-}$ is a bounded quantity, we claim that
\begin{align*}\label{}
  &  \mathop{\sum}_{i=1}^{3} E_{i,\alpha} \leq \mathcal{O}(1)\cdot|\sigma_\alpha|, \ \ \ \ \ \alpha\in NP \\
  &  \mathop{\sum}_{i=1}^{3} E_{i,\alpha}\leq \mathcal{O}(1)\cdot\varepsilon|\alpha|, \ \ \ \ \ \alpha\in S\cup R
\end{align*}
following from the conclusions (3.5)(3.6) in \cite{BressanLiuYang99}.
Recall that $\mathop{\sum}_{\alpha \in NP}|\sigma_\alpha| \leq \varepsilon$. Then it yields
$$ \frac{d}{dt}\Phi(U(t),V(t))=
\mathop{\sum}_{\alpha \in NP}\mathop{\sum}_{i=1}^{3} E_{i,\alpha}+
\mathop{\sum}_{\alpha \in S\cup R}\mathop{\sum}_{i=1}^{3} E_{i,\alpha}
\leq \mathcal{O}(1) \varepsilon.
$$
\hfill{$\square$}

In summary, we draw the following conclusion from Lemmas \ref{Lemma phi decreasing}-
\ref{Lemma phi almost decreasing}.

\begin{Lemma}\label{Lemma stability of approximation}
Let $U^\varepsilon(t),V^\varepsilon(t)$ be $\varepsilon-$approximate solutions of (\ref{ZND}) constructed by fractional step front tracking
method. Choose suitable constants $\epsilon, \kappa_1, \kappa_2, \kappa_3, \kappa_4 $. If
$U^\varepsilon(0), V^\varepsilon(0) \in \mathcal{D}_0,$
 then there exists a constant $C_2$ such that
\begin{align*}
    & \Phi(U^\varepsilon(t),V^\varepsilon(t)) - \Phi(U^\varepsilon(s),V^\varepsilon(s)) \leq C_2 \varepsilon (t-s), \\
    & ||U^\varepsilon(t)-V^\varepsilon(t)||_{L^1} \leq C_1^2||U^\varepsilon(s)-V^\varepsilon(s)||_{L^1}+C_1 C_2\varepsilon (t-s)
\end{align*}
for all $t>s\geq 0,$ where constant $C_1$ is specified in (\ref{equivalence metric}).
\end{Lemma}

\begin{Remark}\label{Remark stability for conservation laws}
Letting $G\equiv \mathbf{0}$ or $Y_0=0$, one can derive the stability estimate for $\mathcal{S}_t^\varepsilon (U)$; that is,
$$||\mathcal{S}_t^\varepsilon (U) - \mathcal{S}_t^\varepsilon (V)||_{L^1}
\leq C_1^2||\mathcal{S}_s^\varepsilon (U)- \mathcal{S}_s^\varepsilon (V) ||_{L^1}+C_1 C_2\varepsilon (t-s)$$
for all $t>s\geq 0$. This coincides with the result of conservation laws in \cite{BressanLiuYang99}.
Besides,  Lipschitz continuity of $\mathcal{S}_t^\varepsilon (U)$ w.r.t time $t$ actually holds.
Therefore we assert that there exists a Lipschitz constant $C_3$ such that
$$||\mathcal{S}_t^\varepsilon (U) - \mathcal{S}_s^\varepsilon (U)||_{L^1}
\leq C_3(t-s)$$
for any $t>s\geq 0$. See \cite{Bressanbook,Dafermosbook}.
\end{Remark}

Lemma \ref{Lemma stability of approximation} and Lemma \ref{Lemma uniqueness} on convergence
give rise to the stability  of entropy solution.


\section{Uniqueness of Limit Solution}\label{section uniqueness}

In the section we attempt to prove that by means of the scheme in Subsection \ref{Def FSFT scheme},
 every sequence of approximate solutions must converge to a unique limit. And then give a rigorous proof
of Theorem 1.

\subsection{Estimates of Operator $\mathcal{P}^\varepsilon_t$}

Note that $\mathcal{S}_t^\varepsilon $ is identity mapping w.r.t. component $Y$.
We claim  the following Lemma developed from \cite{AmadoriGosseGuerra02} still holds
for non-strictly hyperbolic system (\ref{ZND}).
Its proof in detail is referred to Theorem 4.1 in \cite{AmadoriGosseGuerra02}.

\begin{Lemma}\label{Lemma 1}
Assume that $U$ and $V$ are piecewise constant functions in $\mathcal{D}_s$.
$\omega(x) \in L^1(\mathbb{R}^4)$ is a piecewise constant function with small total variation,
and has finitely many discontinuities. Then
$$||\mathcal{S}_t^\varepsilon  (U)-\mathcal{S}_t^\varepsilon  (V)-\omega||_{L^1}\leq
\mathcal{O}(1)|| U- V -\omega||_{L^1} +\mathcal{O}(1) t\{\varepsilon+ T.V.(\omega)\}.$$
\end{Lemma}

Additionally we need some information about source $G$.
Since $G(U)$ in (\ref{ZND})(\ref{Balance Law}) is locally Lipschitz continuous, there exists a constant
$L'>0$ such that
\begin{equation}\label{G Lipschitz continuity}
  \begin{aligned}
     & ||G(U)-G(V)||\leq L' ||U-V||, \\
     & T.V.(G(U)) \leq L'\cdot T.V.(U) , \\
     & ||G(U)||_{L^1} =||G(U)-G(U_\infty)||_{L^1} \leq L'||U-U_\infty||_{L^1}
   \end{aligned}
\end{equation}
for any $U, V \in \mathcal{D}_\infty$.
Next we can calculate the error caused by commutation of $\mathcal{T}_t$
 and $\mathcal{S}_t^\varepsilon$.

\begin{Lemma}\label{Lemma commutation T and S}
(Commutation estimate)
Assume that $U$ is a piecewise constant function in $\mathcal{D}_s$.
Then
$$||\mathcal{S}_t^\varepsilon  \mathcal{T}_t (U)-\mathcal{T}_t\mathcal{S}_t^\varepsilon  (U)||_{L^1}\leq
\mathcal{O}(1)t^2 +\mathcal{O}(1) t \varepsilon$$
holds for $\varepsilon \leq t <\bar{t}$.
\end{Lemma}

{\Proof}
The Lipschitz continuity of $G$ implies
$$T.V.( G(\mathcal{S}_t^\varepsilon  (U)) \leq L'\cdot T.V.(\mathcal{S}_t^\varepsilon(U))
\leq \mathcal{O}(1) \epsilon_\infty.$$
Combining this with Lemma \ref{Lemma 1} and Remark \ref{Remark stability for conservation laws},
we deduce
\begin{align*}
     &||\mathcal{S}_t^\varepsilon \mathcal{T}_t (U)-\mathcal{T}_t\mathcal{S}_t^\varepsilon  (U)||_{L^1}\\
    =& ||\mathcal{S}_t^\varepsilon  (U+tG(U))-\mathcal{S}_t^\varepsilon  (U)-t G(\mathcal{S}_t^\varepsilon  (U))||_{L^1} \\
    \leq &\mathcal{O}(1) || (U+tG(U))-U-t G(\mathcal{S}_t^\varepsilon (U))||_{L^1}
        + \mathcal{O}(1) t \{\varepsilon+ T.V.(t G(\mathcal{S}_t^\varepsilon  (U))\} \\
    \leq &\mathcal{O}(1)t\cdot  ||\mathcal{S}_t^\varepsilon  (U)-U||_{L^1}
         + \mathcal{O}(1) t\cdot  \{\varepsilon+ t\cdot T.V.( G(\mathcal{S}_t^\varepsilon  (U))\} \\
    \leq &\mathcal{O}(1)t^2 + \mathcal{O}(1)t\varepsilon .
\end{align*}
\hfill{$\square$}

\begin{Remark}
Similar property on strictly hyperbolic system was also observed in \cite{ColomboGuerra07}.
Lemma \ref{Lemma commutation T and S} is indeed the basis to establish uniqueness of solutions generated by semigroup $\mathcal{P}$.
It shows that commutation of $\mathcal{S}_\varepsilon^\varepsilon$ and $  \mathcal{T}_\varepsilon$ once gives rise to a small error of $\varepsilon^2$.
Since composition $\mathcal{S}_{k\varepsilon}^\varepsilon(\mathcal{T}_\varepsilon)^k$ can be transformed to $\mathcal{P}_{k\varepsilon}^\varepsilon$ by $k(k-1)/2$ commutations, the total errors accumulated in such process is less than $k^2\varepsilon^2.$ The fact guarantees the conclusion in Lemma \ref{Lemma Pt-ST}.
\end{Remark}

To compute the error of $\mathcal{P}_t^\varepsilon (U)- \mathcal{S}_t^\varepsilon  \mathcal{T}_t (U)$,
 we establish the estimates of operator $(\mathcal{T}_\varepsilon)^k$ in advance.

\begin{Lemma}\label{Lemma T^k estimate}
For any piecewise constant functions $U, V\in \mathcal{D}_s$, the following relations hold.
\begin{align*}
    & ||(\mathcal{T}_\varepsilon)^k (U)-U_\infty||_{L^1}\leq
   \textnormal{e}^{k\varepsilon L'}||U-U_\infty||_{L^1}, \\
   &  ||(\mathcal{T}_\varepsilon)^k (U)-U||_{L^1}\leq
   k\varepsilon L' \textnormal{e}^{k\varepsilon L'}||U-U_\infty||_{L^1}, \\
   &  ||(\mathcal{T}_\varepsilon)^k (U)- (\mathcal{T}_\varepsilon)^k (V) ||_{L^1}\leq
   \textnormal{e}^{k\varepsilon L'}||U-V||_{L^1}.
\end{align*}
\end{Lemma}

{\Proof}
The last inequality is obvious according to Lipschitz continuity of $G$.
It suffices to prove the first two inequalities by induction.
Firstly, for  $k=1$, one deduce from (\ref{G Lipschitz continuity}) that
\begin{align*}
      ||\mathcal{T}_\varepsilon (U)-U_\infty||_{L^1}
    =& ||U-U_\infty+ \varepsilon \{G(U)-G(U_\infty)\}||_{L^1}   \\
    \leq & ||U-U_\infty||_{L^1} + \varepsilon L' ||U-U_\infty||_{L^1} \\
    \leq & \text{e}^{\varepsilon L'} ||U-U_\infty||_{L^1}, \\
      ||\mathcal{T}_\varepsilon (U)-U||_{L^1}
    =& \varepsilon ||G(U)-G(U_\infty)||_{L^1}   \\
    \leq & \varepsilon L' ||U-U_\infty||_{L^1} \\
    \leq & \varepsilon L' \text{e}^{\varepsilon L'} ||U-U_\infty||_{L^1} .
\end{align*}

\noindent Secondly,  for integer $k>1$ , we have
\begin{align*}
      ||(\mathcal{T}_\varepsilon)^k (U)-U_\infty||_{L^1}
    =& ||(\mathcal{T}_\varepsilon)^{k-1} (U)-U_\infty+ \varepsilon
    \{G((\mathcal{T}_\varepsilon)^{k-1} (U))-G(U_\infty)\}||_{L^1}   \\
    \leq & (1+\varepsilon L')||\mathcal{T}_\varepsilon)^{k-1} (U)-U_\infty||_{L^1}  \\
    \leq & \text{e}^{\varepsilon L'} ||\mathcal{T}_\varepsilon)^{k-1} (U)-U_\infty||_{L^1} \\
    \leq & \text{e}^{k \varepsilon L'} ||U-U_\infty||_{L^1}
\end{align*}
and
\begin{align*}
      ||(\mathcal{T}_\varepsilon)^k (U)-U||_{L^1}
    =& ||(\mathcal{T}_\varepsilon)^{k-1} (U)-U+ \varepsilon
    \{G((\mathcal{T}_\varepsilon)^{k-1} (U))-G(U_\infty)\}||_{L^1}   \\
    \leq & ||(\mathcal{T}_\varepsilon)^{k-1} (U)-U ||_{L^1}+
        \varepsilon L' ||(\mathcal{T}_\varepsilon)^{k-1} (U)-U_\infty||_{L^1}  \\
    \leq & (k-1)\varepsilon L'\text{e}^{(k-1)\varepsilon L'} ||U -U_\infty||_{L^1}
        +  \varepsilon L'\text{e}^{(k-1)\varepsilon L'} ||U -U_\infty||_{L^1}    \\
    \leq & k\varepsilon L' \text{e}^{k \varepsilon L'} ||U-U_\infty||_{L^1}
\end{align*}
by induction. \hfill{$\square$}

\begin{Lemma}\label{Lemma Pt-ST}
Assume that $U \in \mathcal{D}_s$ is piecewise constant.
For any time $t<\bar{t}$ and positive $\varepsilon \ll \text{min}\{t,t^2\}$, we have
$$||\mathcal{P}_t^\varepsilon (U)- \mathcal{S}_t^\varepsilon  \mathcal{T}_t (U)||_{L^1}\leq
  \mathcal{O}(1)(1+||U-U_\infty||_{L^1}) t^2 .$$
\end{Lemma}

{\Proof}
Set integer $n=[t/\varepsilon]$ and $k=1,2,\cdots,n.$
Define diameters
\begin{align*}
    &d_k\doteq \mathop{\text{sup}}_{V\in \mathcal{D}_{s+k\varepsilon}}
            ||(\mathcal{T}_\varepsilon)^k  \mathcal{S}_\varepsilon^\varepsilon (V)
            - \mathcal{S}_\varepsilon^\varepsilon  (\mathcal{T}_\varepsilon )^k (V)||_{L^1},\\
    &\bar{d}_k  \doteq \mathop{\text{sup}}_{V\in \mathcal{D}_{s+k\varepsilon}}
             || \mathcal{P}_{k\varepsilon}^\varepsilon (V)
            - \mathcal{S}_\varepsilon^\varepsilon  (\mathcal{T}_\varepsilon )^k (V)||_{L^1}.
\end{align*}

\noindent Applying Lemma \ref{Lemma T^k estimate} and definition of $d_k$, we deduce that
\begin{align*}
    ||(\mathcal{T}_\varepsilon)^k  \mathcal{S}_\varepsilon^\varepsilon (U)
            - \mathcal{S}_\varepsilon^\varepsilon  (\mathcal{T}_\varepsilon )^k (U)||_{L^1}
    \leq & ||(\mathcal{T}_\varepsilon)^{k-1} \mathcal{T}_\varepsilon \mathcal{S}_\varepsilon^\varepsilon (U)
            - (\mathcal{T}_\varepsilon )^{k-1}\mathcal{S}_\varepsilon^\varepsilon  \mathcal{T}_\varepsilon(U)||_{L^1}\\
      &  + || (\mathcal{T}_\varepsilon )^{k-1}\mathcal{S}_\varepsilon^\varepsilon  \mathcal{T}_\varepsilon(U)
        -  \mathcal{S}_\varepsilon^\varepsilon(\mathcal{T}_\varepsilon)^{k-1} \mathcal{T}_\varepsilon (U)||_{L^1} \\
    \leq & \text{e}^{(k-1)\varepsilon L'}|| \mathcal{T}_\varepsilon \mathcal{S}_\varepsilon^\varepsilon (U)
            - \mathcal{S}_\varepsilon ^\varepsilon \mathcal{T}_\varepsilon(U)||_{L^1} \\
           & +|| (\mathcal{T}_\varepsilon )^{k-1}\mathcal{S}_\varepsilon ^\varepsilon \mathcal{T}_\varepsilon(U)
        -  \mathcal{S}_\varepsilon^\varepsilon(\mathcal{T}_\varepsilon)^{k-1} \mathcal{T}_\varepsilon (U)||_{L^1} \\
     \leq & \text{e}^{tL'} d_1 + d_{k-1} .
\end{align*}

\noindent It implies  $d_k  \leq  \text{e}^{tL'} d_1 + d_{k-1}$, furthermore,
\begin{equation}\label{Estimate dk}
 d_k \leq \text{e}^{tL'}k d_1 .
\end{equation}
From stability result in Lemma \ref{Lemma stability of approximation} , Remark \ref{Remark stability for conservation laws} and inequality (\ref{Estimate dk}), we derive
\begin{align*}
    & || \mathcal{P}_{k\varepsilon}^\varepsilon (U) - \mathcal{S}_{k\varepsilon}^\varepsilon  (\mathcal{T}_\varepsilon )^k (U)||_{L^1} \\
   \leq & || \mathcal{P}_{(k-1)\varepsilon}^\varepsilon\mathcal{T}_{\varepsilon}\mathcal{S}_\varepsilon^\varepsilon (U)
     - \mathcal{P}_{(k-1)\varepsilon}^\varepsilon\mathcal{S}_\varepsilon^\varepsilon  \mathcal{T}_\varepsilon(U)||_{L^1}\\
         & \ \ \
     + || \mathcal{P}_{(k-1)\varepsilon}^\varepsilon\mathcal{S}_{\varepsilon}^\varepsilon\mathcal{T}_\varepsilon (U)
     - \mathcal{S}_{(k-1)\varepsilon}^\varepsilon(\mathcal{T}_\varepsilon )^{k-1}\mathcal{S}_\varepsilon^\varepsilon \mathcal{T}_\varepsilon(U)||_{L^1}\\
     & \ \ \
        + || \mathcal{S}_{(k-1)\varepsilon}^\varepsilon(\mathcal{T}_\varepsilon )^{k-1}
        \mathcal{S}_{\varepsilon}^\varepsilon\mathcal{T}_\varepsilon (U)
     - \mathcal{S}_{(k-1)\varepsilon}^\varepsilon\mathcal{S}_\varepsilon^\varepsilon  (\mathcal{T}_\varepsilon )^{k-1} \mathcal{T}_\varepsilon(U)||_{L^1}\\
   \leq & C_1^2|| \mathcal{T}_{\varepsilon}\mathcal{S}_\varepsilon^\varepsilon (U)
     - \mathcal{S}_\varepsilon^\varepsilon  \mathcal{T}_\varepsilon(U)||_{L^1} +  C_1 C_2(k-1)\varepsilon^2 \\
          & \ \ \
     + || \mathcal{P}_{(k-1)\varepsilon}^\varepsilon\mathcal{S}_{\varepsilon}^\varepsilon\mathcal{T}_\varepsilon (U)
     - \mathcal{S}_{(k-1)\varepsilon}^\varepsilon(\mathcal{T}_\varepsilon )^{k-1}\mathcal{S}_\varepsilon^\varepsilon  \mathcal{T}_\varepsilon(U)||_{L^1}  \\
    & \ \ \
     + C_1^2|| (\mathcal{T}_\varepsilon )^{k-1} \mathcal{S}_{\varepsilon}^\varepsilon\mathcal{T}_\varepsilon (U)
     - \mathcal{S}_\varepsilon^\varepsilon  (\mathcal{T}_\varepsilon )^{k-1} \mathcal{T}_\varepsilon(U)||_{L^1}\\
    & \ \ \ + C_1 C_2 (k-1)\varepsilon^2 \\
    \leq & C_1^2 d_1 + \bar{d}_{k-1} + C_1^2 d_{k-1} + 2C_1 C_2 k\varepsilon^2  \\
        \leq & C_1^2(1+\text{e}^{tL'}k) d_1 + \bar{d}_{k-1}  + 2C_1 C_2 k\varepsilon^2 .
\end{align*}

\noindent Hence there holds
$$\bar{d}_{k} \leq C_1^2(1+\text{e}^{tL'}k)  d_1 + \bar{d}_{k-1}+ 2C_1 C_2 k\varepsilon^2,$$
and then
\begin{equation} \label{Estimate d-k}
\begin{aligned}
   \bar{d}_{n} &\leq C_1^2(1+\text{e}^{tL'}n)n  d_1 + 2C_1 C_2 n^2\varepsilon^2 \\
   & = C_1^2\left\{n\varepsilon^2+\text{e}^{tL'}(n\varepsilon)^2\right\}  \frac{d_1}{\varepsilon^2}+
        2C_1 C_2 n^2\varepsilon^2 \\
   & \leq C_1^2(t\varepsilon +\text{e}^{tL'}t^2)  \frac{d_1}{\varepsilon^2}+ 2C_1 C_2 t^2 \\
   & \leq \mathcal{O}(1) t^2,
\end{aligned}
\end{equation}

\noindent since $\varepsilon\leq t^2$ and $d_1 / \varepsilon^2$ is bounded
in light of Lemma \ref{Lemma commutation T and S}.
Observe that Lemma \ref{Lemma T^k estimate} implies
\begin{align*}
    &||(\mathcal{T}_\varepsilon)^k  (U)
            -  \mathcal{T}_{k\varepsilon } (U)||_{L^1} \\
    = & ||\mathcal{T}_\varepsilon (\mathcal{T}_\varepsilon)^{k-1}   (U)
            - \mathcal{T}_{(k-1)\varepsilon} (U) - \varepsilon G(U)||_{L^1} \\
    = &  ||(\mathcal{T}_\varepsilon)^{k-1} (U)
            - \mathcal{T}_{(k-1)\varepsilon} (U)
            + \varepsilon G((\mathcal{T}_\varepsilon)^{k-1} (U))- \varepsilon G(U)||_{L^1} \\
     \leq & ||(\mathcal{T}_\varepsilon)^{k-1} (U)- \mathcal{T}_{(k-1)\varepsilon} (U) ||_{L^1}
            + (k-1)\varepsilon^2 L'^2 \text{e}^{(k-1)\varepsilon L'} || U-  U_\infty||_{L^1}.
\end{align*}
Hence
$$||(\mathcal{T}_\varepsilon)^k  (U) -  \mathcal{T}_{k\varepsilon } (U)||_{L^1}
\leq k^2\varepsilon^2 L'^2 \text{e}^{k\varepsilon L'} || U-  U_\infty||_{L^1} ,$$
moreover,
\begin{equation}\label{Estimate Tn-Tt}
\begin{aligned}
    &||(\mathcal{T}_\varepsilon)^n  (U)  -  \mathcal{T}_{t} (U)||_{L^1} \\
  \leq & ||(\mathcal{T}_\varepsilon)^n  (U)  -  \mathcal{T}_{n\varepsilon} (U)||_{L^1}
        + ||\mathcal{T}_{t} (U) -  \mathcal{T}_{n\varepsilon} (U)||_{L^1} \\
  \leq & n^2\varepsilon^2 L'^2 \text{e}^{n\varepsilon L'} || U-  U_\infty||_{L^1}
        +(t-n\varepsilon) L' || U-  U_\infty||_{L^1}   \\
  \leq & 2 t^2 L'^2 \text{e}^{t L'}|| U-  U_\infty||_{L^1}
\end{aligned}
\end{equation}
provided $L'>1$.
Consequently, we see that
\begin{align*}
    ||\mathcal{P}_t^\varepsilon (U)- \mathcal{S}_t^\varepsilon  \mathcal{T}_t (U)||_{L^1}
    \leq & ||\mathcal{P}_t^\varepsilon (U)-   \mathcal{P}_{n\varepsilon}^\varepsilon (U)||_{L^1}
        + ||\mathcal{P}_{n\varepsilon}^\varepsilon (U)- \mathcal{S}_{n\varepsilon}^\varepsilon  \mathcal{T}_{n\varepsilon} (U)||_{L^1}\\
     &\ \ \   + ||\mathcal{S}_{n\varepsilon}^\varepsilon  \mathcal{T}_{n\varepsilon} (U) - \mathcal{S}_{t}^\varepsilon  \mathcal{T}_{t} (U)||_{L^1} \\
    \leq & \mathcal{O}(1)\varepsilon^2 + \bar{d}_n + \mathcal{O}(1)||\mathcal{T}_{n\varepsilon} (U) - \mathcal{T}_{t} (U)||_{L^1} \\
      &\ \ \       + \mathcal{O}(1)(t-n\varepsilon) \varepsilon \\
    \leq & \mathcal{O}(1)\varepsilon^2 +\mathcal{O}(1) (1+ || U-  U_\infty||_{L^1})t^2  \\
    \leq & \mathcal{O}(1) (1+ || U-  U_\infty||_{L^1}) t^2
\end{align*}
by Lemma \ref{Lemma stability of approximation}, Remark \ref{Remark stability for conservation laws}
and estimates (\ref{Estimate d-k})(\ref{Estimate Tn-Tt}).
\hfill{$\square$}

In the end we state an approximate semigroup property of $\mathcal{P}_t^\varepsilon$.

\begin{Lemma}\label{Lemma Semigroup property}
(Semigroup estimate)
For any $t_1, t_2 >0$ and $U\in \mathcal{D}_s$, the inequality
$$||\mathcal{P}_{t_2}^\varepsilon \mathcal{P}_{t_1}^\varepsilon (U)- \mathcal{P}_{t_1 +t_2}^\varepsilon (U)||_{L^1}\leq
  \mathcal{O}(1)(1+\varepsilon +t_1 + t_2+||U-U_\infty||_{L^1}) \varepsilon$$
holds if $0<\varepsilon \ll\bar{t}$.
\end{Lemma}

{\Proof}
Fix $\varepsilon \ll\bar{t}$. Then there exists two integers
$k$ and $K$ such that
$$(k-1)\varepsilon < t_1 \leq k\varepsilon, \ \ \ \ K\varepsilon \leq t_1 +t_2 < (K+1)\varepsilon.$$
They yield
$$\left[\frac{t_2}{\varepsilon}\right] = K-k \ \ \text{or} \ \ K-k+1.$$
Since the arguments of these two cases are analogous, we only discuss the case  $[t_2/\varepsilon] = K-k $.

Let $V\doteq \mathcal{P}_{t_1}^\varepsilon (U) ,  \ V_\infty:= \mathcal{P}_{t_1}^\varepsilon (U_\infty) $. Then
Lemma \ref{Lemma stability of approximation} and Remark \ref{Remark stability for conservation laws}
give that
\begin{align*}
  &  ||\mathcal{P}_{t_2}^\varepsilon \mathcal{P}_{t_1}^\varepsilon (U)- \mathcal{P}_{t_1 +t_2}^\varepsilon (U)||_{L^1}  \\
\leq & ||\mathcal{S}_{t_1+t_2-K\varepsilon}^\varepsilon \mathcal{S}_{k\varepsilon- t_1}^\varepsilon
        \mathcal{P}_{(K-k)\varepsilon}^\varepsilon (V) -
        \mathcal{S}_{t_1+t_2-K\varepsilon}^\varepsilon  \mathcal{P}_{(K-k)\varepsilon}^\varepsilon \mathcal{T}_\varepsilon \mathcal{S}_{k\varepsilon- t_1}^\varepsilon (V)
    ||_{L^1}  \\
\leq & \mathcal{O}(1)|| \mathcal{S}_{k\varepsilon- t_1}^\varepsilon \mathcal{P}_{(K-k)\varepsilon}^\varepsilon (V)
        -\mathcal{P}_{(K-k)\varepsilon}^\varepsilon \mathcal{T}_\varepsilon \mathcal{S}_{k\varepsilon- t_1}^\varepsilon (V)
    ||_{L^1} \\
     & \ \ +  \mathcal{O}(1) (t_1+t_2-K\varepsilon)\varepsilon\\
\leq & \mathcal{O}(1)|| \mathcal{P}_{(K-k)\varepsilon}^\varepsilon \mathcal{T}_\varepsilon \mathcal{S}_{k\varepsilon- t_1}^\varepsilon (V)
    -  \mathcal{P}_{(K-k)\varepsilon}^\varepsilon (V)||_{L^1} \\
    & \ \
    +  \mathcal{O}(1)|| \mathcal{S}_{k\varepsilon- t_1}^\varepsilon \mathcal{P}_{(K-k)\varepsilon}^\varepsilon (V)
    -  \mathcal{P}_{(K-k)\varepsilon}^\varepsilon (V)||_{L^1} + \mathcal{O}(1) \varepsilon^2 \\
\leq & \mathcal{O}(1)||  \mathcal{T}_\varepsilon \mathcal{S}_{k\varepsilon- t_1}^\varepsilon (V) -  V||_{L^1}
    + \mathcal{O}(1) (K-k)\varepsilon^2 \\
    & \ \ +  \mathcal{O}(1) (k\varepsilon - t_1) + \mathcal{O}(1) \varepsilon^2 \\
\leq & \mathcal{O}(1)\{|| \mathcal{S}_{k\varepsilon- t_1}^\varepsilon (V) -  V||_{L^1}
    +  \varepsilon || G(\mathcal{S}_{k\varepsilon- t_1}^\varepsilon (V)) -  G(\mathcal{S}_{k\varepsilon- t_1}^\varepsilon (V_\infty))||_{L^1} \} \\
    & \ \ \    + \mathcal{O}(1) (t_2\varepsilon + \varepsilon+ \varepsilon^2)  \\
\leq &  \mathcal{O}(1) ||V-V_\infty||_{L^1}  \varepsilon
       + \mathcal{O}(1)(k\varepsilon - t_1)\varepsilon
    +    \mathcal{O}(1)( t_2\varepsilon + \varepsilon+ \varepsilon^2)  \\
\leq & \mathcal{O}(1)||U-U_\infty||_{L^1}  \varepsilon + \mathcal{O}(1) (t_1+t_2 +1+\varepsilon)\varepsilon .
\end{align*}
We complete the proof.\hfill{$\square$}

\subsection{Convergence of $\varepsilon-$approximate Solutions}

First, we need a local estimate of different trajectories $\mathcal{P}_t^\mu$
and  $\mathcal{P}_t^\nu$ which originate from the same state $U$.

\begin{Lemma}\label{lemma local error Pt}
Assume that $U \in \mathcal{D}_s$ is piecewise constant . Then the inequality
$$||\mathcal{P}_t^\mu (U)- \mathcal{P}_t^\nu (U)||_{L^1}\leq
  \mathcal{O}(1)(1+||U-U_\infty||_{L^1}) t^2$$
holds for any positive $\mu,\ \nu \ll \textrm{min}\{t,t^2\}$, and time $t<\bar{t}$.
\end{Lemma}

{\Proof}
Note that for any $V\in \mathcal{D}_\infty$,
$$||\mathcal{S}_t^\mu (V)- \mathcal{S}_t^\nu  (V)||_{L^1}
\leq \mathcal{O}(1)\text{max}\{\mu, \nu\}t$$
holds (cf. \cite{Bressanbook}) .
By Lemma \ref{Lemma Pt-ST} we deduce
\begin{align*}
    &||\mathcal{P}_t^\mu (U)- \mathcal{P}_t^\nu (U)||_{L^1} \\
\leq & ||\mathcal{P}_t^\mu (U)- \mathcal{S}_t^\mu \mathcal{T}_t(U)||_{L^1}
       + ||\mathcal{P}_t^\nu (U)- \mathcal{S}_t^\nu \mathcal{T}_t(U)||_{L^1}
       +  ||\mathcal{S}_t^\mu \mathcal{T}_t (U)- \mathcal{S}_t^\nu \mathcal{T}_t(U)||_{L^1} \\
\leq & \mathcal{O}(1)(1+||U-U_\infty||_{L^1}) t^2 + \mathcal{O}(1)\text{max}\{\mu, \nu\}t \\
\leq & \mathcal{O}(1)(1+||U-U_\infty||_{L^1}) t^2  .
\end{align*}
\hfill{$\square$}

Second, we can generalize the result in Lemma \ref{lemma local error Pt} to that
of trajectories with distinct origins.
The following Lemma exhibits  that $\{\mathcal{P}_t^\varepsilon (U^\varepsilon)\}$
is a Cauchy sequence  w.r.t. $\varepsilon$ in $L^1$-topology.

\begin{Lemma}\label{Lemma uniqueness}
Assume that $U^\mu_0, U^\nu_0$ are piecewise constant functions in $ \mathcal{D}_0$ such that
$$ ||U^\mu_0- U_0||_{L^1} < \mu, \ \ \ \ ||U^\nu_0- U_0||_{L^1} < \nu ,$$
$$T.V.(U^\mu_0) \leq T.V.(U_0), \ \ \ \ T.V.(U^\nu_0) \leq T.V.(U_0)$$
with some $U_0 \in \mathcal{D}_0$. Then for every $t\geq 0$,
\begin{equation}\label{convergence P}
  ||\mathcal{P}_t^\mu (U^\mu_0)-\mathcal{P}_t^\nu (U^\nu_0)||_{L^1} \rightarrow 0
\end{equation}
as $\mu, \nu \rightarrow 0$.
\end{Lemma}

{\Proof}
Divide the interval $[0,t]$ equally into $m$ parts such that $\Delta t= t/m$ is small.
Let the points of such division  satisfy
$$0=t_0 < t_1 < \cdots < t_m = t.$$
If $\mu, \nu < \text{min}\{\bar{t},(\Delta t) ^2\}$,  it follows from Lemma \ref{Lemma Semigroup property}
that
\begin{equation}\label{Estimate 4}
 \begin{aligned}
    &  ||\mathcal{P}_t^\mu (U^\mu_0)-\mathcal{P}_t^\nu (U^\mu_0)||_{L^1}  \\
\leq & \mathop{\sum}_{i=1}^m ||\mathcal{P}_{t-t_i}^\mu \mathcal{P}_{t_i}^\nu(U^\mu_0)
        -\mathcal{P}_{t-t_{i-1}}^\mu \mathcal{P}_{t_{i-1}}^\nu(U^\mu_0)||_{L^1} \\
\leq & \mathop{\sum}_{i=1}^m ||\mathcal{P}_{t-t_i}^\mu \mathcal{P}_{t_i}^\nu(U^\mu_0)
        - \mathcal{P}_{t-t_{i}}^\mu\mathcal{P}_{t_i-t_{i-1}}^\mu \mathcal{P}_{t_{i-1}}^\nu(U^\mu_0)||_{L^1}  \\
    & \ \ \ +  \mathcal{O}(1) \mathop{\sum}_{i=1}^m
      (1+\mu+ t-t_{i-1} + ||\mathcal{P}_{t_{i-1}}^\nu(U^\mu_0) - U_\infty||_{L^1})\mu \\
=: & \  \Sigma_1 + \Sigma_2
\end{aligned}
\end{equation}
where $\Sigma_1$ (resp. $\Sigma_2$) denotes the first (resp. the second) sum of the
 third part in estimate (\ref{Estimate 4}).
Thus using Lemmas \ref{Lemma Semigroup property},\ref{lemma local error Pt},\ref{Lemma stability of approximation}
and Remark \ref{Remark stability for conservation laws},
we obtain

\begin{equation}\label{Estimate sigma1}
\begin{aligned}
    \Sigma_1
    & \leq \mathcal{O}(1)  \mathop{\sum}_{i=1}^m || \mathcal{P}_{t_i}^\nu(U^\mu_0)
        -\mathcal{P}_{t_i-t_{i-1}}^\mu \mathcal{P}_{t_{i-1}}^\nu(U^\mu_0)||_{L^1}
         + \mathcal{O}(1) \mathop{\sum}_{i=1}^m (t-t_i)\mu \\
    & \leq \mathcal{O}(1)  \mathop{\sum}_{i=1}^m || \mathcal{P}_{t_i-t_{i-1}}^\nu\mathcal{P}_{t_{i-1}}^\nu(U^\mu_0)
        -\mathcal{P}_{t_i-t_{i-1}}^\mu \mathcal{P}_{t_{i-1}}^\nu(U^\mu_0)||_{L^1}  \\
    & \ \ \      + \mathcal{O}(1)  \mathop{\sum}_{i=1}^m (1+ \nu + t_i   +||U^\mu_0 - U_\infty||_{L^1}) \nu
         +\mathcal{O}(1) \mathop{\sum}_{i=1}^m (t-t_i)\mu \\
    & \leq \mathcal{O}(1)  \mathop{\sum}_{i=1}^m (1+ ||\mathcal{P}_{t_{i-1}}^\nu(U^\mu_0) - U_\infty||_{L^1} )(t_i-t_{i-1})^2 \\
    & \ \ \      + \mathcal{O}(1)  (1+ \nu + t  +||U^\mu_0 - U_\infty||_{L^1})m \nu
         +\mathcal{O}(1) tm\mu \\
    &  \leq \mathcal{O}(1)  \mathop{\sum}_{i=1}^m (1+ ||\mathcal{P}_{t_{i-1}}^\nu(U^\mu_0) -\mathcal{P}_{t_{i-1}}^\nu(U_\infty)||_{L^1}+
            ||\mathcal{S}_{t_{i-1}}^\nu(U_\infty)- U_\infty||_{L^1} )(\Delta t)^2 \\
    & \ \ \      + \mathcal{O}(1)  (1+ \nu + t +||U^\mu_0 - U_\infty||_{L^1}  )m \nu
         +\mathcal{O}(1) tm\mu \\
    &  \leq \mathcal{O}(1)  \mathop{\sum}_{i=1}^m (1+ \nu t_{i-1} +t_{i-1}+
            ||U^\mu_0- U_\infty||_{L^1} )(\Delta t)^2 \\
    & \ \ \      + \mathcal{O}(1)  (1+ \nu +  t +||U^\mu_0 - U_\infty||_{L^1}  )m \nu
         +\mathcal{O}(1) tm\mu \\
    &  \leq \mathcal{O}(1)   (1+ \mu + \nu t +t+
            ||U_0- U_\infty||_{L^1} )t\Delta t \\
    & \ \ \      + \mathcal{O}(1)  (1+\mu+ \nu + t   +||U_0 - U_\infty||_{L^1})m \nu
         +\mathcal{O}(1) tm\mu ,\\
\end{aligned}
\end{equation}

\begin{equation}\label{Estimate sigma2}
\begin{aligned}
    \Sigma_2
    & \leq \mathcal{O}(1)\mathop{\sum}_{i=1}^m  (1+\mu+ t-t_{i-1} + ||\mathcal{P}_{t_{i-1}}^\nu(U^\mu_0) -
       \mathcal{P}_{t_{i-1}}^\nu(U_\infty)||_{L^1}  \\
    &  \ \ \ \ \ \ \ \ \ \ \ \ \ \ \ \ \ \
         + ||\mathcal{S}_{t_{i-1}}^\nu(U_\infty)- U_\infty||_{L^1})\mu   \\
    & \leq \mathcal{O}(1) \mathop{\sum}_{i=1}^m (1+\mu+ t-t_{i-1} + t_{i-1}+ \nu t_{i-1}  +||U^\mu_0 - U_\infty||_{L^1})\mu   \\
    & \leq \mathcal{O}(1)  (1+\mu  + t + \nu t  +||U_0 - U_\infty||_{L^1})m \mu .
\end{aligned}
\end{equation}

\noindent Estimates (\ref{Estimate 4})-(\ref{Estimate sigma2}) yield
\begin{align*}
   &  ||\mathcal{P}_t^\mu (U^\mu_0)-\mathcal{P}_t^\nu (U^\mu_0)||_{L^1}\\
 \leq & \mathcal{O}(1)   (1+ \mu +\nu +t+ \nu  t+  ||U_0- U_\infty||_{L^1} )\{t\Delta t +m (\mu+\nu)\}.
\end{align*}

\noindent As a result,
\begin{align*}
    &    ||\mathcal{P}_t^\mu (U^\mu_0)-\mathcal{P}_t^\nu (U^\nu_0)||_{L^1} \\
\leq & ||\mathcal{P}_t^\mu (U^\mu_0)-\mathcal{P}_t^\nu (U^\mu_0)||_{L^1}
        + ||\mathcal{P}_t^\nu (U^\mu_0)-\mathcal{P}_t^\nu (U^\nu_0)||_{L^1} \\
\leq &  \mathcal{O}(1)    (1+ \mu +\nu +t+ \nu t+  ||U_0- U_\infty||_{L^1} )\{t\Delta t +m (\mu+\nu)\} \\
    & \ \ \    + \mathcal{O}(1) || U^\mu_0 -U^\nu_0||_{L^1} + \mathcal{O}(1)t\nu \\
\leq &  \mathcal{O}(1)    (1+ \mu +\nu +t+ \nu t+  ||U_0- U_\infty||_{L^1} )\{t\Delta t +m (\mu+\nu)\} .
\end{align*}

\noindent Passing $\mu, \nu \rightarrow 0$, we have
$$\mathop{\overline{\text{lim}}}_{\mu, \nu \rightarrow 0} ||\mathcal{P}_t^\mu (U^\mu_0)-\mathcal{P}_t^\nu (U^\nu_0)||_{L^1}
\leq  \mathcal{O}(1)   (1  +t+  ||U_0- U_\infty||_{L^1} )t\Delta t ,$$
which implies (\ref{convergence P}) due to the arbitrariness of $\Delta t$.
\hfill{$\square$}

In the light of previous preparation, we continue to prove Theorem 1.\\

\noindent \textbf{Proof of Theorem 1}

Assume that $\epsilon$ is small enough, and $U_0 (x) \in \mathcal{D}_0$.
Then according to Lemma \ref{Lemma uniqueness}, every sequence of $\varepsilon-$approximate
solution $\mathcal{P}_t^\varepsilon (U_0^\varepsilon)$ to Cauchy problem (\ref{ZND})
(\ref{initial data}) converges to
a unique limit solution $U(t,x) \doteq \mathcal{P}_t (U_0)$.
It was proved in \cite{Hu18} $U(t,x)$ is an entropy solution to problem
(\ref{ZND})(\ref{initial data}).
Thus it suffices to verify that  the mapping $\mathcal{P} :\ [0,+\infty)\times\mathcal{D}_\infty
\rightarrow \mathcal{D}_\infty$
is a Lipschitz semigroup.

Semigroup property (\ref{Semigroup P}) is deduced from uniqueness and
Lemmas \ref{Lemma Semigroup property} and \ref{Lemma uniqueness}.
Now we only check the Lipschitz continuity of $\mathcal{P}$.
At every time $t_k =k \varepsilon$, one see that
\begin{align*}
  ||\mathcal{P}_{t_k}^{\varepsilon}(U_0^{\varepsilon}) - \mathcal{P}_{t_k-}^{\varepsilon}(U_0^{\varepsilon})||_{L^1}
   &
  = \mathcal{O}(1) \int_{-\infty}^{+\infty} Y(t_k -,x) \phi(T(t_k-,x))\varepsilon dx  \\
  & \leq  \mathcal{O}(1) \text{e}^{-\phi k\varepsilon} \varepsilon \int_{-\infty}^{+\infty}
   Y(0,x)dx
\end{align*}
according to definition (\ref{U(t_k,x)}) and decay rate in (\ref{Decay of Y(t)}).
By finite propagation speed of waves, we assert that
$$ ||\mathcal{P}_{t'}^{\varepsilon}(U_0^{\varepsilon}) - \mathcal{P}_{s'}^{\varepsilon}(U_0^{\varepsilon})||_{L^1} \leq \mathcal{O}(1)
 T.V.(\mathcal{P}_{t_{k-1}}^{\varepsilon}(U_0^{\varepsilon})) \cdot|t'-s'|$$
for any $t',s' \in [t_{k-1}, t_k).$
Suppose  $t>s\geq 0$. Choose $\varepsilon$ small enough, then have
\begin{align*}
  ||\mathcal{P}_{t}^{\varepsilon}(U_0^{\varepsilon}) - \mathcal{P}_{s}^{\varepsilon}(U_0^{\varepsilon})||_{L^1}
   &  \leq  \mathcal{O}(1) \left(t-s+  \mathop{\sum}\limits_{s\leq k\varepsilon \leq t }
    \text{e}^{-\phi k\varepsilon}  \varepsilon  \right) \\
  & \leq \mathcal{O}(1) \left(t-s+ \int_s^t \text{e}^{-\phi\tau} d\tau\right) \\
  & = \mathcal{O}(1)  |t-s|.
\end{align*}

\noindent As $\varepsilon \rightarrow 0$, there holds
$$||\mathcal{P}_t (U_0)- \mathcal{P}_s (U_0)||_{L^1}  \leq \mathcal{O}(1)  |t-s|.$$
Remember that Lemma \ref{Lemma stability of approximation} implies for any $U_0, V_0 \in \mathcal{D}_0$,
$$||\mathcal{P}_s (U_0)- \mathcal{P}_s (V_0)||_{L^1}  \leq \mathcal{O}(1) ||U_0- V_0||_{L^1} $$
if passing to limits. Finally preceding two inequalities yield (\ref{Lipschitz P}).
\hfill{$\square$}

\section{Local Characteristics of Entropy Weak Solutions}

This section is concerned with the local geometric features of entropy solutions.
$U_0^\varepsilon$ is as usual  regarded as the $\varepsilon-$approximation of initial data $U_0\in \mathcal{D}_0$.
Based on Theorem 1 in this paper and the convergence of
$\mathcal{S}_t^\varepsilon$ in \cite{Bressanbook},
we claim that in $L^1$-topology sense,
\begin{align*}
     U_0^\varepsilon & \rightarrow U_0,  \\
     \mathcal{P}_t^\varepsilon (U_0^\varepsilon)  & \rightarrow \mathcal{P}_t (U_0) , \\
     \mathcal{S}_t^\varepsilon \mathcal{T}_t(U_0^\varepsilon) & \rightarrow \mathcal{S}_t\mathcal{T}_t(U_0), \\
     \mathcal{T}_t\mathcal{S}_t^\varepsilon (U_0^\varepsilon)&  \rightarrow \mathcal{T}_t\mathcal{S}_t (U_0),
\end{align*}
as $\varepsilon \rightarrow 0$.
 Hence by Lemmas \ref{Lemma commutation T and S} and \ref{Lemma Pt-ST}, we draw the following conclusion.

\begin{Lemma}
Suppose $U(s,x)=\mathcal{P}_s (U_0)$ with $U_0\in \mathcal{D}_0$. Then
\begin{align*}
   & ||\mathcal{P}_\theta (U)-\mathcal{S}_\theta \mathcal{T}_\theta(U)||_{L^1} \leq
 \mathcal{O}(1)(1+||U-U_\infty||_{L^1}) \theta^2 , \\
   & ||\mathcal{P}_\theta (U)-\mathcal{T}_\theta \mathcal{S}_\theta(U)||_{L^1} \leq
 \mathcal{O}(1)(1+||U-U_\infty||_{L^1}) \theta^2, \\
   & ||\mathcal{S}_\theta \mathcal{T}_\theta(U)-\mathcal{T}_\theta \mathcal{S}_\theta(U)||_{L^1} \leq
 \mathcal{O}(1) \theta^2
\end{align*}
hold for  $\theta \in (0, \bar{t})$.
\end{Lemma}

We will utilize this Lemma to analyze the relations between balance laws (\ref{ZND})
and conservation laws (resp. transport system) .
Now let us focus on a small neighborhood of point $(s,\xi)$. Introduce Riemann problem of conservation laws,
\begin{equation}\label{conservation law}
    \begin{aligned}
        & w_t+ F(w)_x =0, \\
        & w(0,x) = \left\{
                \begin{aligned}
                    & U(s, \xi+)\ \ \ \text{if} \ x<\xi, \\
                    & U(s, \xi-)\ \ \ \text{if} \ x>\xi.
                \end{aligned}\right.
    \end{aligned}
\end{equation}

\noindent and Cauchy problem of transport equations.
\begin{equation}\label{transport eq}
    \begin{aligned}
        & w_t+ dF(U(s, \xi))w_x = G(U(s,x)), \\
        & w(0,x) = U(s,x).
    \end{aligned}
\end{equation}

\noindent Then regard $U_\xi^C$ as the entropy solution to (\ref{conservation law}),
and $U_\xi^T$ as the solution to (\ref{transport eq}).
One can figure out that
\begin{align*}
    & U_\xi^T (t,x)= U_\xi^1 + U_\xi^2, \\
    & U_\xi^1 (t,x)= \mathop{\sum}_{i=1}^3  (l_i^\xi \cdot U(s,x - \lambda_i^\xi t))  r_i^\xi,  \\
    & U_\xi^2 (t,x)= \mathop{\sum}_{i=1}^3  \int_0^t (l_i^\xi \cdot G(U(s,x - \lambda_i^\xi \tau)) r_i^\xi d\tau,
\end{align*}
where $l_i^\xi$ (resp. $r_i^\xi$) denotes the left (resp. right) eigenvector of $dF(U(s, \xi))$.
Now we proceed to prove Theorem 2.\\

\noindent \textbf{Proof of Theorem 2}

Suppose $U(t,x)= \mathcal{P}_t (U_0)$ with $U_0 \in \mathcal{D}_0$. It is an entropy solution to (\ref{ZND})(\ref{initial data}).
Fix $s\geq0$ and $\xi \in \mathbb{R}$.
Then we compute the limit (\ref{Local character 1}).
Observe that
\begin{align*}
    & \frac{1}{\theta}\int_{\xi-\theta \hat{\lambda}}^{\xi+\theta \hat{\lambda}}
    ||U(s+\theta, x)-U_\xi^C (\theta, x)||dx  \\
 \leq  &  \frac{1}{\theta}\int_{\xi-\theta \hat{\lambda}}^{\xi+\theta \hat{\lambda}}
            ||\mathcal{P}_\theta (U(s, x))-S_\theta (U(s, x)) - \theta G(U(s, x))||dx \\
       &  + \frac{1}{\theta}\int_{\xi-\theta \hat{\lambda}}^{\xi+\theta \hat{\lambda}}
            ||S_\theta ( U(s, x)) - U_\xi^C (\theta, x)||dx   +
            \frac{1}{\theta}\int_{\xi-\theta \hat{\lambda}}^{\xi+\theta \hat{\lambda}}
            ||\theta G(U(s, x))||dx  \\
 \leq  &  \frac{1}{\theta} ||\mathcal{P}_\theta (U(s, x))-S_\theta (U(s, x)) - \theta G(U(s, x))||_{L^1}\\
       &  + \frac{1}{\theta}\int_{\xi-\theta \hat{\lambda}}^{\xi+\theta \hat{\lambda}}
            ||S_\theta ( U(s, x)) - U_\xi^C (\theta, x)||dx   +
            \int_{\xi-\theta \hat{\lambda}}^{\xi+\theta \hat{\lambda}} ||G(U(s, x))||dx  .
\end{align*}

\noindent Since
\begin{align*}
    ||\mathcal{P}_\theta (U)-\mathcal{T}_t \mathcal{S}_\theta(U)||_{L^1} \leq
 \mathcal{O}(1)(1+||U-U_\infty||_{L^1}) \theta^2
\end{align*}
for  $\theta \in (0, \bar{t})$, we have
\begin{align*}
&  \mathop{\text{lim}}_{\theta \rightarrow 0^+}
   \frac{1}{\theta} ||\mathcal{P}_\theta (U(s, x))-S_\theta (U(s, x)) - \theta G(U(s, x))||_{L^1} \\
= &  \mathop{\text{lim}}_{\theta \rightarrow 0^+}
\frac{1}{\theta} ||\mathcal{P}_\theta (U(s, x))-\mathcal{T}_\theta \mathcal{S}_ \theta(U(s, x)) ||_{L^1} \\
=& 0
\end{align*}

\noindent Additionally
$$ \mathop{\text{lim}}_{\theta \rightarrow 0^+}
\frac{1}{\theta}\int_{\xi-\theta \hat{\lambda}}^{\xi+\theta \hat{\lambda}}
 ||S_\theta ( U(s, x)) - U_\xi^C (\theta, x)||dx =0  $$
is proved by Bressan \cite{Bressanbook}(see also \cite{Dafermosbook});
$$ \mathop{\text{lim}}_{\theta \rightarrow 0^+}
 \int_{\xi-\theta \hat{\lambda}}^{\xi+\theta \hat{\lambda}} ||G(U(s, x))||dx  =0$$
is deduced from absolute continuity of Lebesgue integral.
Therefore conclusion (\ref{Local character 1}) holds.

Next we calculate the limit (\ref{Local character 2}). For every $a<\xi<b$,
\begin{align*}
&    \frac{1}{\theta}\int_{a+\theta \hat{\lambda}}^{b-\theta \hat{\lambda}}
        ||U(s+\theta,x)-U_\xi^T(\theta,x)||dx    \\
\leq  & \frac{1}{\theta} \int_{a+\theta \hat{\lambda}}^{b-\theta \hat{\lambda}}
        ||\mathcal{P}_\theta (U(s))(x)-\mathcal{S}_\theta (U(s, x)) - \theta G(U(s, x))||dx  \\
     &+ \frac{1}{\theta} \int_{a+\theta \hat{\lambda}}^{b-\theta \hat{\lambda}}
        ||\mathcal{S}_\theta (U(s, x)) - U_\xi^1 (\theta,x)||dx  \\
     &+ \frac{1}{\theta} \int_{a+\theta \hat{\lambda}}^{b-\theta \hat{\lambda}}
        ||\theta G(U(s, x)) - U_\xi^2 (\theta,x)||dx  \\
\leq   & \frac{1}{\theta}
        ||\mathcal{P}_\theta (U(s))(x)-\mathcal{T}_\theta \mathcal{S}_\theta (U(s, x))||_{L^1}  + \mathcal{O}(1) \{T.V.(U(s); (a,b))\}^2  \\
     &+ \frac{1}{\theta} \int_{a+\theta \hat{\lambda}}^{b-\theta \hat{\lambda}}
        ||\theta G(U(s, x)) - U_\xi^2 (\theta,x)||dx.
\end{align*}
The first term  of the sum above tends to 0 as $\theta\rightarrow 0^+$.
The second term was derived by \cite{AmadoriGosseGuerra02,Bressanbook} etc..
So it suffices to compute the third term .
In fact
\begin{align*}
  &  \frac{1}{\theta} \int_{a+\theta \hat{\lambda}}^{b-\theta \hat{\lambda}}
        ||\theta G(U(s, x)) - U_\xi^2 (\theta,x)||dx    \\
\leq & \frac{1}{\theta} \int_{a+\theta \hat{\lambda}}^{b-\theta \hat{\lambda}}
        \left\|\int_0^\theta \mathop{\sum}_{i=1}^3 l_i^\xi \cdot G(U(s,x )) r_i^\xi d\tau -
        \int_0^\theta \mathop{\sum}_{i=1}^3 l_i^\xi \cdot G(U(s,x - \lambda_i^\xi \tau)) r_i^\xi d\tau
        \right\|dx  \\
\leq & \mathcal{O}(1) \frac{1}{\theta}\mathop{\sum}_{i=1}^3
        \int_0^\theta \int_{a+\theta \hat{\lambda}}^{b-\theta \hat{\lambda}}
        \left\|   G(U(s,x ))  -   G(U(s,x - \lambda_i^\xi \tau)) \right\|dx   d\tau  \\
\leq & \mathcal{O}(1) \frac{1}{\theta}\mathop{\sum}_{i=1}^3
        \int_0^\theta  T.V. (G(U(s));(a,b)) \cdot  \lambda_i^\xi \tau  d\tau  \\
\leq & \mathcal{O}(1) \theta \cdot T.V. (G(U(s));(a,b)).
\end{align*}

\noindent Consequently we have
$$ \mathop{\text{lim}}_{\theta \rightarrow 0^+}
\frac{1}{\theta} \int_{a+\theta \hat{\lambda}}^{b-\theta \hat{\lambda}}
        ||\theta G(U(s, x)) - U_\xi^2 (\theta,x)||dx =0.$$
Then (\ref{Local character 2}) holds for some constant $C_0$.

Conversely, suppose mapping $U:\ [0,+\infty)
\rightarrow \mathcal{D}_\infty $  is  Lipschitz continuous in $L^1-$topology and satisfies
(\ref{Local character 1})(\ref{Local character 2}) for $s \in [0,+\infty)\setminus\mathcal{N}$ and $\xi\in\mathbb{R}$,
where $\mathcal{N}$ is a set  with measure 0. Given $s \in [0,+\infty)\setminus\mathcal{N}$ and $\varepsilon>0$,
we can cover the interval $[a,b]$ with finitely many  intervals $[a_i,b_i]$ such that
\begin{itemize}
  \item $\bigcup\limits_i [a_i,b_i]=[a,b]$  and $(a_i,b_i)\bigcap(a_j,b_j)=\emptyset $ with $i\neq j$;
  \item $T.V.(U(s);(a_i,b_i))<\varepsilon$ ;
  \item $\mathop{\overline{\textnormal{lim}}}_{\theta \rightarrow 0^+}
        \frac{1}{\theta}\int_{a_i+\theta \hat{\lambda}}^{b_i-\theta \hat{\lambda}}
        ||U(s+\theta,x)-U_{\zeta_i}^T(\theta,x)||dx  \leq
            C_0 \{T.V.(U(s);(a_i,b_i))\}^2$
        where $\zeta_i \doteq (a_i+b_i)/2$.
\end{itemize}

\noindent Set $\xi_i=a_i$. Then construct a function
\begin{equation}\label{}
  U^*(t,x)\doteq \left\{
        \begin{aligned}
           & U_{\xi_i}^C(t,x) &\ \ \ \ &\text{if}\ \ |x-\xi_i|\leq t\hat{\lambda} , \\
           & U_{\zeta_i}^T(t,x) &\ \ \ \ &\text{if}\ \ a_i+t\hat{\lambda}<x<b_i-t\hat{\lambda} ,
        \end{aligned}
        \right.
\end{equation}
for $t$ small enough. Hence
\begin{align*}
  &\mathop{\overline{\textnormal{lim}}}_{\theta \rightarrow 0^+}
\frac{1}{\theta}\int_{a+\theta \hat{\lambda}}^{b-\theta \hat{\lambda}}
||U(s+\theta)(x)-U^*(\theta,x)||dx  \\
\leq  & \sum\limits_i \mathop{\textnormal{lim}}_{\theta \rightarrow 0^+}
        \frac{1}{\theta}\int_{\xi_i-\theta \hat{\lambda}}^{\xi_i+\theta \hat{\lambda}}
        ||U(s+\theta)(x)-U_{\xi_i}^C(\theta,x)||dx  \\
      & \ \ + \sum\limits_i \mathop{\overline{\textnormal{lim}}}_{\theta \rightarrow 0^+}
        \frac{1}{\theta}\int_{a_i+\theta \hat{\lambda}}^{b_i-\theta \hat{\lambda}}
        ||U(s+\theta)(x)-U_{\zeta_i}^T(\theta,x)||dx   \\
\leq  & \sum\limits_i
         C_0 \{T.V.(U(s);(a_i,b_i))\}^2        \\
\leq  &  C_0  T.V.(U(s);(a,b)) \varepsilon .
\end{align*}

\noindent Similarly we have
$$\mathop{\overline{\textnormal{lim}}}_{\theta \rightarrow 0^+}
\frac{1}{\theta}\int_{a+\theta \hat{\lambda}}^{b-\theta \hat{\lambda}}
||\mathcal{P}_{s+\theta} (U_0)-U^*(\theta,x)||dx
\leq    C_0 T.V.(\mathcal{P}_{s} (U_0);(a,b))\varepsilon   . $$
Therefore,
\begin{align*}
   & \mathop{\overline{\textnormal{lim}}}_{\theta \rightarrow 0^+}
\frac{1}{\theta}\int_{a+\theta \hat{\lambda}}^{b-\theta \hat{\lambda}}
||U(s+\theta)(x) - \mathcal{P}_{s+\theta} (U_0)||dx  \\
\leq  & C_0  \{T.V.(U(s);(a,b)) +  T.V.(\mathcal{P}_{s} (U_0);(a,b))\} \varepsilon.
\end{align*}
We find $U(s )= \mathcal{P}_{s} (U_0)$ according to the arbitrariness of $(a,b)$ and $\varepsilon$ in the last estimate.
\hfill{$\square$}



\bibliography{mybibfile}

\begin{thebibliography}{10}
\expandafter\ifx\csname url\endcsname\relax
  \def\url#1{\texttt{#1}}\fi
\expandafter\ifx\csname urlprefix\endcsname\relax\def\urlprefix{URL }\fi
\expandafter\ifx\csname href\endcsname\relax
  \def\href#1#2{#2} \def\path#1{#1}\fi

\bibitem{DafermosHsiao82}
C.~Dafermos, L.~Hsiao, Hyperbolic systems of balance laws with inhomogeneity
  and dissipation, Indiana Univ. Math. J. 31~(4) (1982) 471--491.

\bibitem{AmadoriGosseGuerra02}
D.~Amadori, L.~Gosse, G.~Guerra, Global bv entropy solutions and uniqueness for
  hyperbolic systems of balance laws, Arch. Rational Mech. Anal. 162~(4) (2002)
  327--366.

\bibitem{AmadoriGuerra02}
D.~Amadori, G.~Guerra, Uniqueness and continuous dependence for systems of
  balance laws with dissipation, Nonlinear Analysis 49 (2002) 987--1014.

\bibitem{AmadoriGuerra99}
D.~Amadori, G.~Guerra, Global weak solutions for systems of balance laws, Appl.
  Math. Lett. 12~(6) (1999) 123--127.

\bibitem{Christoforou06}
C.~Christoforou, Hyperbolic systems of balance laws via vanishing viscosity,
  Journal of Differential Equations 221~(2) (2006) 470--541.

\bibitem{ColomboGuerra07}
R.~Colombo, G.Guerra, Hyperbolic balance laws with a non local source,
  Communications in Partial Differential Equations 32~(12) (2007) 1917--1939.

\bibitem{ChenWagner03}
G.~Chen, D.~Wagner, Global entropy solutions to exothermically reacting,
  compressible euler equations, J. Differential Equations 191~(2) (2003)
  277--322.

\bibitem{Hu18}
K.~Hu, Existence of global bv solutions to a model of reacting euler fluid with
  variable thermodynamics parameters, J. Math.Anal.Appl. 457 (2018) 890--921.

\bibitem{Chen92}
G.~Chen, Global solutions to the compressible navier-stokes equations for a
  reacting mixture, SIAM J. Math. Anal. 23~(3) (1992) 609--634.

\bibitem{ChenHoffTrivisa03}
G.~Chen, D.~Hoff, K.~Trivisa, Global solutions to a model for exothermically
  reacting, compressible flows with large discontinuous initial data, Arch.
  Rational Mech. Anal. 166~(4) (2003) 321--358.

\bibitem{Wang04}
D.~Wang, Global solution for the mixture of real compressible reacting flows in
  combustion, Commun. Pure Appl. Anal 3~(4) (2004) 775--790.

\bibitem{DonatelliTrivisa06}
D.Donatelli, K.Trivisa, On the motion of a viscous compressible
  radiative-reacting gas, Commun. Math. Phys. 265~(2) (2006) 463--491.

\bibitem{BressanLiuYang99}
A.~Bressan, T.P.Liu, T.Yang, L$^1$ stability estimates for $n\times n$
  conservation laws, Arch. Rational Mech. Anal. 149 (1999) 1--22.

\bibitem{LiuYang99}
T.P.Liu, T.Yang, L$^1$ stability for $2\times2$ systems of hyperbolic
  conservation laws, J. Amer. Math. Soc. 12 (1999) 729--774.

\bibitem{Bressanbook}
A.~Bressan, Hyperbolic systems of conservation laws: the one-dimensional cauchy
  problem, Oxford University Press.

\bibitem{Dafermosbook}
C.~Dafermos, Hyperbolic conservation laws in continuum physics, Springer.

\end{thebibliography}

\end{document}